\documentclass{amsart}
\usepackage{amsmath}
\usepackage{amssymb}
\usepackage{amsthm}

\DeclareMathOperator{\sech}{sech}

\DeclareMathOperator{\tr}{tr}
\DeclareMathOperator{\Vol}{Vol}
\DeclareMathOperator{\dvol}{dvol}

\DeclareMathOperator{\Ric}{Ric}

\DeclareMathOperator{\connectsum}{\#}

\DeclareMathOperator{\Met}{Met}

\newcommand{\og}{\overline{g}}

\newcommand{\ophi}{\overline{\phi}}

\newcommand{\lp}{\langle}
\newcommand{\rp}{\rangle}
\newcommand{\lv}{\lvert}
\newcommand{\rv}{\rvert}

\newcommand{\auxiliarymetric}{auxiliary metric}
\newcommand{\auxiliarymanifold}{auxiliary manifold}

\newcommand{\charconstant}{characteristic constant}

\newcommand{\mB}{\mathcal{B}}
\newcommand{\mC}{\mathcal{C}}
\newcommand{\mD}{\mathcal{D}}

\newcommand{\mF}{\mathcal{F}}

\newcommand{\mW}{\mathcal{W}}

\newcommand{\kM}{\mathfrak{M}}

\newcommand{\kX}{\mathfrak{X}}

\newcommand{\bC}{\mathbb{C}}

\newcommand{\bN}{\mathbb{N}}

\newcommand{\bR}{\mathbb{R}}

\newcommand{\bZ}{\mathbb{Z}}

\newcommand{\comment}[1]{}

\newtheorem{thm}{Theorem}[section]
\newtheorem{prop}[thm]{Proposition}
\newtheorem{lem}[thm]{Lemma}
\newtheorem{cor}[thm]{Corollary}

\theoremstyle{definition}
\newtheorem{defn}[thm]{Definition}

\newtheorem{example}[thm]{Example}

\theoremstyle{remark}
\newtheorem{remark}[thm]{Remark}

\numberwithin{equation}{section}

\begin{document}

\title[SMMS and Quasi-Einstein Metrics]{Smooth metric measure spaces and quasi-Einstein metrics}
\author{Jeffrey S. Case}
\thanks{Partially supported by NSF-DMS Grant No.\ 1004394}
\address{Department of Mathematics \\ Princeton University \\ Princeton, NJ 08544}
\email{jscase@math.princeton.edu}
\date{}
\keywords{Einstein, conformally Einstein, quasi-Einstein, gradient Ricci soliton, smooth metric measure space}
\subjclass[2000]{Primary 53Cxx; Secondary 53A30,53C25}
\begin{abstract}
Smooth metric measure spaces have been studied from the two different perspectives of Bakry-\'Emery and Chang-Gursky-Yang, both of which are closely related to work of Perelman on the Ricci flow.  These perspectives include a generalization of the Ricci curvature and the associated quasi-Einstein metrics, which include Einstein metrics, conformally Einstein metrics, gradient Ricci solitons, and static metrics.  In this article, we describe a natural perspective on smooth metric measure spaces from the point of view of conformal geometry and show how it unites these earlier perspectives within a unified framework.  We offer many results and interpretations which illustrate the unifying nature of this perspective, including a natural variational characterization of quasi-Einstein metrics as well as some interesting families of examples of such metrics.
\end{abstract}
\maketitle

\section{Introduction}
\label{sec:intro}

Over the past quarter century, smooth metric measure spaces have begun to attract a lot of attention in Riemannian geometry.  These spaces are Riemannian manifolds $(M^n,g)$ equipped with a smooth measure $e^{-\phi}\dvol_g$.  The natural question to ask about these spaces is how the change of measure should affect their study as geometric objects.

Among the most basic, and obvious, changes one must make is to find a suitable generalization of the Ricci curvature.  This is because the Ricci curvature controls the rate of volume growth, so if one changes the meaning of the volume, of course one must also change the meaning of the Ricci curvature.  Essentially motivated by the Bochner (in)equality, Bakry and \'Emery~\cite{BakryEmery1985} proposed that the natural definition should be
\[ \Ric_\phi^m := \Ric + \nabla^2\phi - \frac{1}{m}d\phi\otimes d\phi, \]
where the constant $m\in[0,\infty]$ specifies the ``dimension'' of the measure $e^{-\phi}\dvol_g$.  Indeed, one can show that if $\Ric_\phi^m\geq Kg$, then the volume of a geodesic ball of radius $R$ with respect to $e^{-\phi}\dvol_g$ is less than that of a ball of the same radius in the $(m+n)$-dimensional spaceform of constant sectional curvature $\frac{K}{m+n-1}$.  The argument is based on the Bochner inequality, and as such one can extend a great many results from comparison geometry in the Riemannian setting to the setting of smooth metric measure spaces; for details, see~\cite{Wei_Wylie} and the references therein.

Another important application of the Bakry-\'Emery Ricci tensor is as a description of quasi-Einstein metrics; i.e.\ Riemannian manifolds $(M^n,g)$ such that there exists a function $\phi$ and constants $m\in[0,\infty]$, $\lambda\in\bR$ such that $\Ric_\phi^m=\lambda g$.  Special cases include Einstein metrics, static metrics, gradient Ricci solitons, and more generally, the bases of Einstein warped products~\cite{Besse,Case_Shu_Wei,HePetersenWylie2010} for various positive values of $m$, and also conformally Einstein metrics and conformally Einstein products~\cite{Chen2009} when $m$ is negative.  Due to the importance of all of these metrics, and the fact that they are related by varying the constant $m$, it is natural to regard the notion of a quasi-Einstein metric as ``interpolating'' between these different metrics.

In the special case $m=\infty$, smooth metric measure spaces played an important role in Perelman's work~\cite{Perelman1} on the Ricci flow.  One of Perelman's first key observations was that the Ricci flow can be realized as a gradient flow in the space of smooth metric measure spaces.  Indeed, Perelman found a generalization of the scalar curvature for smooth metric measure spaces with the property that gradient Ricci solitons $\Ric_\phi^\infty=\lambda g$ are critical points of the natural total scalar curvature functional $\mF$.  Using the fact that $\Ric_\phi^\infty$ differs from the Ricci curvature by a Lie derivative, it then follows that the gradient flow for this functional is the Ricci flow, modulo a measure-dependent diffeomorphism.

Recently, a different perspective on smooth metric measure spaces was taken up by Chang, Gursky and Yang~\cite{CGY0}.  They were interested in determining whether one can meaningfully discuss conformal invariants of smooth metric measure spaces, starting with the question of whether there were conformally invariant notions of the scalar and Ricci curvatures.  For trivial reasons, this turns out to be the case.  However, this observation has the somewhat surprising corollary that both $\Ric_\phi^\infty$ and Perelman's $\mF$-functional can be regarded as the $m\to\infty$ limits of the corresponding conformally invariant Ricci curvature and total scalar curvatures.

These different works raise a number of interesting questions about smooth metric measure spaces, all of which we aim to address in this work.  The main question is whether the perspective centered around the study of the Bakry-\'Emery Ricci tensor and the perspective centered around finding conformal invariants can be brought together as a part of one larger framework for studying smooth metric measure spaces.  As we will see, the answer is yes, and it is based upon two observations.  First, there is a natural notion of a conformal transformation of a smooth metric measure space.  This notion is implicit in the work of Chang, Gursky and Yang~\cite{CGY0}, but it seems not to have been made explicit anywhere in the literature.  Second, this notion of a conformal transformation leads to a notion of a ``duality'' between smooth metric measure spaces with two specifically related values of $m$.  In particular, this duality relates the perspectives of Bakry and \'Emery~\cite{BakryEmery1985} and of Chang, Gursky and Yang~\cite{CGY0}.

Having unified these two perspectives, one wants to explore further the consequences of the ideas of conformal transformations and the interpolating properties of $m$, which is the second aim of this article.  To address the first point, we give explicit formulae for how the Bakry-\'Emery Ricci tensor and the natural scalar curvature of a smooth metric measure space change under a conformal transformation, and use it to extend some estimates from~\cite{Case_Shu_Wei} to the case of quasi-Einstein metrics with $m<1-n$.  This relationship will also manifest itself in our definition of the energy functional of a smooth metric measure space, which we will observe in~\cite{Case2010b} to also provide a way to interpolate between the Yamabe functional and Perelman's $\mW$-functional.  To address the second point, we discuss how many well-known facts for gradient Ricci solitons can be realized as limits as $m\to\infty$ of facts known for Einstein metrics on products.  In particular, we include a discussion of some well-known examples of Einstein metrics and gradient Ricci solitons, and show how they can be regarded as one-parameter families of quasi-Einstein metrics parameterized by $m\in(1,\infty]$.

In the sequel~\cite{Case2010b}, we expand upon both of these points and prove a precompactness theorem for the space of quasi-Einstein metrics on a Riemannian manifold with $m\in(1,\infty]$.  Our result allows $m$ to vary, showing that the examples discussed in this article are in some sense typical of quasi-Einstein metrics.  A key ingredient in our proof will be to introduce and understand some basic consequences of the generalization of the Yamabe constant to the setting of smooth metric measure spaces, expanding upon the observed relation between the Yamabe constant and Perelman's $\lambda$-entropy made by Chang, Gursky and Yang~\cite{CGY0}.  In particular, we will also establish a similar understanding of Perelman's $\nu$-entropy, and use it to prove a $\kappa$-noncollapsing-type result for compact quasi-Einstein metrics.

This article is organized as follows.

In Section~\ref{sec:background}, we provide more details about how smooth metric measure spaces are studied in the spirits of Bakry-\'Emery, of Perelman, and of Chang-Gursky-Yang.  In particular, we will highlight the key ways in which these perspectives use the notion of a smooth metric measure space so that we might incorporate them into our more general framework.

In Section~\ref{sec:smooth}, we describe our proposal for the study of quasi-Einstein metrics.  First, we give a precise definition of a smooth metric measure space which specifically emphasizes the role of the parameter $m$.  We also define what we mean by a conformal transformation of a smooth metric measure space.  We then use these definitions to describe our program for the study of smooth metric measure spaces by finding natural invariants which are the ``weighted'' analogues of geometric notions important in Riemannian and conformal geometry.

In Section~\ref{sec:quasi_einstein}, we give the definitions of the Bakry-\'Emery Ricci tensor and the weighted scalar curvature, and use them to define and discuss quasi-Einstein metrics.  By considering how these notions of curvature transform under conformal changes of a smooth metric measure space, we arrive at the aforementioned notion of duality.  We will also observe that the weighted scalar curvature leads to the natural variational characterization of quasi-Einstein metrics, and recall and establish some simple estimates for compact quasi-Einstein metrics.

In Section~\ref{sec:examples}, we conclude by discussing some interesting families of quasi-Einstein metrics.  The first nontrivial examples we discuss are the so-called elliptic Gaussians, which are in many ways the model spaces for quasi-Einstein smooth metrics measure spaces.  We then discuss rotationally symmetric quasi-Einstein smooth metric measure spaces satisfying $\Ric_\phi^m=0$, which are well-understood in both the context of warped products and of gradient Ricci solitons.  In particular, this discussion allows us to realize Hamilton's cigar soliton~\cite{Hamilton1986b} as the limit as $m\to\infty$ of the two-dimensional bases of the $(m+1)$-dimensional Riemannian Schwarzschild metric, and the Bryant soliton as the limit of examples of Einstein metrics constructed by B\"ohm~\cite{Bohm1999}.  We also discuss the Einstein warped products constructed by L\"u, Page and Pope~\cite{LPP} and show how, in the four-dimensional case, they limit to the Cao-Koiso~\cite{Cao1996,Koiso1990} shrinking gradient Ricci soliton as $m\to\infty$.  Finally, we discuss two simple product constructions for smooth metric measure spaces which can be used to create new quasi-Einstein metrics.

\section*{Acknowledgments}
This article is based in part upon my Ph.D.\ dissertation~\cite{Case_dissertation}, completed under the supervision of Xianzhe Dai, who started this work by asking me how the perspectives of Bakry-\'Emery and Chang-Gursky-Yang are related.  Additionally, many of the ideas and their presentation have benefited greatly from conversations with Robert Bartnik, Rod Gover, Chenxu He, Pengzi Miao, Peter Petersen, Yujen Shu, Guofang Wei and William Wylie.  Finally, I wish to thank Michele Rimoldi for comments which helped improve an early version of this article.

\section{Background}
\label{sec:background}

To motivate our perspective, let us first discuss in more detail the existing perspectives on smooth metric measure spaces existing in the literature.

\subsection{The Point of View of Bakry and \'Emery}
\label{sec:warped}

The starting point of Bakry and \'Emery's work on smooth metric measure spaces $(M^n,g,e^{-\phi}\dvol_g)$ is the $\phi$-Laplacian $\Delta_\phi=\Delta-\nabla\phi\cdot$, which is formally self-adjoint with respect to the measure $e^{-\phi}\dvol_g$.  Given an effective dimension $N=m+n$, $m\in[0,\infty]$, one defines the Bakry-\'Emery Ricci tensor as the curvature term of the Bochner (in)equality
\begin{align*}
\frac{1}{2}\Delta_\phi\lv\nabla w\rv^2 & = \lv\nabla^2w\rv^2 + \frac{1}{m}(\lp\nabla\phi,\nabla w\rp)^2 + \lp\nabla w,\nabla\Delta_\phi w\rp + \Ric_\phi^m(\nabla w,\nabla w) \\
& \geq \frac{1}{m+n}(\Delta_\phi w)^2 + \lp\nabla w,\nabla\Delta_\phi w\rp + \Ric_\phi^m(\nabla w,\nabla w) .
\end{align*}
As the Bochner formula underpins many results of comparison geometry, this definition allows one to discuss comparison results for the Bakry-\'Emery Ricci tensor for arbitrary values of $m\geq0$ (cf.\ \cite{BakryQian2000,Qian1997,Wei_Wylie}).

On the other hand, given two Riemannian manifolds $(M,g)$ and $(F^m,h)$ and some smooth function $f\in C^\infty(M)$, it is easy to see that the sequence of warped products
\begin{equation}
\label{eqn:bg_wp}
\left( M\times F, g\oplus (\varepsilon e^{-f/m})^2 h\right)
\end{equation}
converges in the measured Gromov-Hausdorff topology to the smooth metric measure space $(M,g,e^{-\phi}\dvol_g)$ as $\varepsilon\to0$, where the measure arises as the renormalized Riemannian measure on the warped product (cf.\ \cite{CheegerColding1997}).  In this setting, it is straightforward to check that for any vector field $X$ which is tangent to the base $M$ and for any fixed $\varepsilon>0$, the Ricci curvature of the warped product manifold~\eqref{eqn:bg_wp} in the direction $X$ is given by $\Ric_\phi^m(X,X)$.  In particular, this gives a natural ``extrinsic'' way to extend results from comparison geometry to smooth metric measure spaces (cf.\ \cite{Lott2003}).  That is to say, we can either use the Bakry-\'Emery Ricci tensor to study smooth metric measure spaces intrinsically, or we can introduce the \auxiliarymanifold\ $M\times F$ via~\eqref{eqn:bg_wp} to study them extrinsically.

\subsection{The Point of View of Perelman}
\label{sec:grs}

To see how the Ricci flow arises as a gradient flow, let $(M,g,e^{-\phi}\dvol_g)$ be a compact smooth metric measure space such that $e^{-\phi}\dvol_g$ is a probability measure.  Define the functional $\mF$ on such pairs $(g,e^{-\phi}\dvol_g)$ by
\[ \mF(g,e^{-\phi}\dvol_g) = \int_M (R+\lv\nabla\phi\rv^2) e^{-\phi}\dvol_g . \]
Fixing the measure $e^{-\phi}\dvol_g$, a straightforward computation (cf.\ \cite{Perelman1}) yields that the gradient flow of $\mF$ is given by
\begin{equation}
\label{eqn:gradient_flow}
\begin{cases}
          \frac{\partial g}{\partial t} & = -2(\Ric + \nabla^2\phi) \\
          \frac{\partial \phi}{\partial t} & = -(R+\Delta \phi) ,
\end{cases}
\end{equation}
where $\phi(t)$ is defined by $e^{-\phi}\dvol_g=e^{-\phi(t)}\dvol_{g(t)}$.  Note in particular the appearance of the Bakry-\'Emery Ricci tensor $\Ric_\phi:=\Ric_\phi^\infty$.

In order to recover the ordinary Ricci flow, let $(g(t),\phi(t))$ be a solution of \eqref{eqn:gradient_flow} and let $\psi(t)$ be the family of diffeomorphisms generated by $\nabla\phi(t)$.  If $\og(t)=(\psi(t))^\ast g(t)$ and $\ophi(t) = (\psi(t))^\ast\phi(t)$, one finds that $(\og(t),\ophi(t))$ solves
\[ \begin{cases}
          \frac{\partial\og}{\partial t} & = -2\Ric \\
          \frac{\partial\ophi}{\partial t} & = -(R+\Delta \ophi - \lv\nabla \ophi\rv^2) .
   \end{cases} \]
In particular, $\og(t)$ is a solution of the Ricci flow.

As was pointed out by Perelman, different choices for the measure $e^{-\phi}\dvol_g$ lead to the same flow, up to diffeomorphism; the difference lies in the solvability of the system~\eqref{eqn:gradient_flow}.  In this context, we may thus think of a choice of measure as being equivalent to a choice of gauge.

In allowing one to view the Ricci flow as a gradient flow, the notion of a smooth metric measure space also allows one to find many monotone quantities along the flow, among them Perelman's entropies.  These entropies are closely related to sharp logarithmic Sobolev inequalities, and play a key role in establishing $\kappa$-noncollapsing results for the Ricci flow~\cite{Perelman1}.

Finally, we note that Perelman also demonstrated how his heuristic argument establishing the monotonicity of the reduced length and volume can be phrased in terms of the measure $e^{-\phi}\dvol_g$.  These two notions were necessary to establish certain comparison-type theorems for smooth metric measure spaces which cannot be achieved using the Bakry-\'Emery Ricci tensor (cf.\ \cite{Wei_Wylie}).  Unfortunately, at this stage we do not know if the reduced length and volume admit natural analogues for smooth metric measure spaces with finite $m$.
\subsection{The Point of View of Chang, Gursky and Yang}
\label{sec:conformal}

In~\cite{CGY0}, Chang, Gursky and Yang asked if, given a smooth metric measure space $(M^n,g,d\omega)$ with the measure $d\omega=u^{-n}\dvol_g$ regarded as being fixed, one can define conformally invariant analogues of the Ricci and scalar curvatures.  In a certain sense, the answer is trivially ``yes:'' Fixing the measure $d\omega$ imposes the requirement that a function $w$ act as a conformal transformation of $(M^n,g,d\omega)$ by
\[ \left(M^n,g,d\omega\right) \mapsto \left( M^n,w^{-2}g,d\omega \right) . \]
In particular, since $d\omega=\dvol_{u^{-2}g}$, one sees that the measure defines a canonical metric in the conformal class $[g]$, and one thus defines the curvatures
\begin{align*}
\Ric_\omega & := \Ric_{u^{-2}g} \\
R_\omega & = u^{-2}R_{u^{-2}g} .
\end{align*}
By construction, these are conformally invariant tensors which generalize the Ricci and scalar curvature, corresponding to the case $d\omega=\dvol_g$.

From this basic idea, Chang, Gursky and Yang proceed to make two quite interesting observations.  First, similar to the discussion in Section~\ref{sec:warped}, one can introduce a dimensional parameter $m$ and define the ``conformally invariant'' Ricci curvature $\Ric_\omega^m$ of $(M^n,g,d\omega)$ as the horizontal component of $\Ric_{\omega\times\dvol_h}$ on the smooth metric measure space $(M^n\times T^m,g\oplus h,d\omega\times\dvol_h)$.  Of course, $\Ric_\omega^m$ is not a conformal invariant of $[g]$, but rather of the conformal class $[g\oplus h]$ on $M\times T^m$.  Writing $d\omega=u^{-m-n}\dvol_g$ and $u^{m+n-2}=e^f$, this definition is such that
\[ \Ric_\omega^m = \Ric + \nabla^2 f + \frac{1}{m+n-2}df\otimes df + \frac{1}{m+n-2}\left(\Delta f - \lv\nabla f\rv^2\right) g . \]
In particular, as $m\to\infty$, the tensor $\Ric_\omega^m$ converges to the Bakry-\'Emery Ricci tensor $\Ric_f^\infty$ of the smooth metric measure space $(M^n,g,e^{-f}\dvol_g)$.  They define the scalar curvature $R_\omega^m$ similarly, and observe that it converges as $m\to\infty$ to the scalar curvature of the smooth metric measure space $(M^n,g,e^{-f}\dvol_g)$ introduced by Perelman~\cite{Perelman1}.

Having a natural ``conformally invariant'' scalar curvature on a smooth metric measure space, Chang, Gursky and Yang then posed the question of studying the ``Yamabe functional''
\[ \mF(g,d\omega) = \int_M u^2R_\omega^m d\omega = \int_M \left(R + \frac{m+n-1}{m+n-2}\lv\nabla f\rv^2\right) e^{-f}\dvol_g, \]
where we have written $d\omega=u^{-m-n}\dvol_g=e^{-\frac{m+n}{m+n-2}f}\dvol_g$ as above.  When $m=0$, this is of course the usual Yamabe functional, while when $m=\infty$ it is Perelman's $\mF$-functional, thereby establishing an interesting, if not completely understood, link between these two important functionals.

Second, using a result of Moser~\cite{Moser1965}, they showed that the usual Yamabe problem is equivalent to the problem of minimizing $\mF$ over a diffeomorphism class; i.e.\ the Yamabe constant $\sigma\left([g]\right)$ is equivalently written
\[ \sigma\left([g]\right) = \inf_{\psi\in\mD} \mF\left(\psi^\ast g,\omega\right) , \]
where $\mD$ is the diffeomorphism group of $M$, thereby giving an alternative perspective on Perelman's use of the measure to fix a gauge.           
\section{Smooth Metric Measure Spaces}
\label{sec:smooth}

The perspective on smooth metric measure spaces we take throughout this article is based upon the following definition.

\begin{defn}
A \emph{smooth metric measure space (SMMS)} is a four-tuple $(M^n,g,v^m\dvol_g,m)$ of a Riemannian manifold $(M^n,g)$ together with its canonical Riemannian volume element $\dvol_g$, a positive function $v\in C^\infty(M)$, and a dimensional parameter $m\in\bR\cup\{\pm\infty\}$.

Equivalently, a SMMS is a four-tuple $(M^n,g,e^{-\phi}\dvol_g,m)$ of a Riemannian manifold $(M^n,g)$, a function $\phi\in C^\infty(M)$, and a dimensional parameter $m\in\bR\cup\{\pm\infty\}$.
\end{defn}

That these definitions are equivalent is seen by defining $v^m=e^{-\phi}$ for $m\not\in\{\pm\infty\}$.  From now on, we shall always reserve the symbols $v$ and $\phi$ to denote the measure of a SMMS, and they will always be related by $v^m=e^{-\phi}$.  Furthermore, we shall frequently simply denote a SMMS by the triple $(M^n,g,v^m\dvol)$, with the dimensional parameter $m$ understood as being determined by the exponent of $v$; when $\lv m\rv=\infty$, we make sense of this by writing $v^{\pm\infty}=e^{-\phi}$.  Finally, when $m=0$ our convention is such that the measure is always $v^0\dvol_g=\dvol_g$; i.e.\ we are simply studying a Riemannian manifold.

The role of the dimensional parameter $m$ is to specify that the measure $v^m\dvol$ should be regarded as an $(m+n)$-dimensional measure.  This is made precisely through the definition of natural geometric quantities associated to a SMMS, as we will describe below.  One way to motivate this perspective is by considering special sequences of collapsing Riemannian manifolds.

\begin{example}
\label{ex:cc}
Let $(M^n,g,v^m\dvol)$ be a SMMS with $m\in\bN$, and let $(F^m,h)$ be any Riemannian manifold.  For each $\varepsilon>0$, consider the warped product
\begin{equation}
\label{eqn:cc_wp}
\left( M^n\times F^m, g\oplus (\varepsilon v)^2h \right) .
\end{equation}
It is straightforward to check that
\[ \left( M^n\times F^m, g\oplus (\varepsilon v)^2h, \overline{\dvol_{g\oplus v^2h}}, (p,q) \right) \to (M^n, g, v^m\dvol, p) \]
in the pointed measured Gromov-Hausdorff topology as $\varepsilon\to 0$, where $\overline{\dvol_{g\oplus v^2h}}$ denotes the renormalized volume element of~\eqref{eqn:cc_wp}; see~\cite{CheegerColding1997} for details.
\end{example}

In this way, one can take the perspective that a SMMS is what remains of a collapsing sequence of warped products $M\times_v F^m$ as in~\eqref{eqn:cc_wp}, except that the notion of a SMMS allows the dimension of the fiber to be any real number, or even infinite.  As we will see, many geometric notions associated to a SMMS we introduce can be realized by considering the warped product $M\times_v F^m$ with this freedom of specifying $m$.  To make this precise, it is convenient to introduce the following terminology:

\begin{defn}
Let $(M^n,g,v^m\dvol_g,m)$ be a SMMS with $m\in\bN$ and fix a Riemannian manifold $(F^m,h)$.  The \emph{\auxiliarymanifold} is the warped product manifold
\[ \left( M^n\times F^m, g\oplus v^2 h\right) . \]
\end{defn}

Our objective is to study SMMS as geometric objects.  To accomplish this, we will often need to determine what are the ``weighted'' analogues of geometric notions associated to a Riemannian manifold.  After the \emph{total weighted volume}
\[ \Vol_\phi(M) := \int_M v^m\dvol_g, \]
the most natural such notion is the weighted divergence.

\begin{defn}
Let $(M^n,g,v^m\dvol_g)$ be a SMMS and let $(V,h_V)$ and $(W,h_W)$ be vector bundles with inner product over $M$, and let $\lp\cdot,\cdot\rp_V$ and $\lp\cdot,\cdot\rp_W$ be the corresponding inner products
\begin{align*}
\lp\zeta_1,\zeta_2\rp_V & = \int_M h_V(\zeta_1,\zeta_2) v^m\dvol_g \\
\lp\xi_1,\xi_2\rp_W & = \int_M h_W(\xi_1,\xi_2) v^m\dvol_g
\end{align*}
determined by the measure $e^{-\phi}\dvol_g$, where $\zeta_i\in\Gamma(V)$, $\xi_i\in\Gamma(W)$, $i=1,2$ are sections of $V$ and $W$, respectively.  The \emph{weighted divergence $\delta_\phi\colon\Gamma(W)\to\Gamma(V)$} of an operator $D\colon\Gamma(V)\to\Gamma(W)$ is the (negative of the) formal adjoint of $D$ with respect to the inner products $\lp\cdot,\cdot\rp_V$ and $\lp\cdot,\cdot\rp_W$; i.e.\ for all $\zeta\in\Gamma(V)$ and $\xi\in\Gamma(W)$, at least one of which is compactly supported in $M$,
\[ \lp D(\zeta), \xi\rp_W = -\lp \zeta, \delta_\phi\xi\rp_V . \]

The \emph{weighted Laplacian $\Delta_\phi\colon C^\infty(M)\to C^\infty(M)$} is the operator $\Delta_\phi=\delta_\phi d$.
\end{defn}

We opt to use the terminology ``weighted divergence,'' rather than the more appropriate ``weighted adjoint,'' because we shall only consider here the cases where $\delta_\phi$ is defined in terms of the exterior derivative or the Levi-Civita connection, where the terminology divergence is standard.

A straightforward consequence of the definition of the weighted divergence is the following useful formula for computing it in terms of the usual divergence.

\begin{lem}
\label{lem:weighted_divergence}
Let $(M^n,g,v^m\dvol_g)$ be a SMMS.  The weighted divergence $\delta_\phi$ is related to the usual divergence $\delta$ by
\[ \delta_\phi = e^\phi\circ\delta\circ e^{-\phi}, \]
where $e^\phi$ and $e^{-\phi}$ are regarded as multiplication operators.  In particular, we have the formulae
\begin{align*}
\delta_\phi\omega & = \delta\omega - \imath_{\nabla\phi} \omega \\
\Delta_\phi w & = \Delta w - \lp\nabla\phi,\nabla w\rp
\end{align*}
for all $\omega\in\Lambda^kT^\ast M$ and all $w\in C^\infty(M)$.
\end{lem}

As stated in the introduction, one of the two main purposes of this article is to describe how ideas from conformal geometry naturally enter into the study of SMMS.  In order to accomplish this, we must first define what it means for two SMMS to be conformally equivalent.

\begin{defn}
Two SMMS $(M^n,g,v^m\dvol_g)$ and $(M^n,\hat g,\hat v^m\dvol_{\hat g})$ are said to be \emph{(pointwise) conformally equivalent} if there is a positive function $u\in C^\infty(M)$ such that
\begin{equation}
\label{eqn:scms}
\left(M^n,\hat g,\hat v^m\dvol_{\hat g}\right) = \left( M^n,u^{-2}g, u^{-m-n}v^m\dvol_g\right) .
\end{equation}
\end{defn}

In particular, observe that the measure transforms as would the Riemannian volume element of an $(m+n)$-dimensional manifold.

To make sense of this definition in the limit $\lv m\rv=\infty$, one defines the function $f$ by $u^{m+n-2}=e^f$.  In this way, one sees that a ``conformal transformation'' of a SMMS with $\lv m\rv=\infty$ is simply a change of measure.  It is in this way which we will understand many definitions made in the context of the Ricci flow as natural analogues of definitions made in conformal geometry.  Throughout this article, whenever $u$ is thought of as a conformal factor, we will always associated to it the function $f$ defined by $u^{m+n-2}=e^f$.

As in any subject, one wants to study natural invariants on SMMS.  As geometric objects, the way to make precise what is a natural invariant is by using the action of the diffeomorphism group.

\begin{defn}
Let $(M^n,g,v^m\dvol_g)$ be a SMMS.  An object $T=T[g,v^m\dvol_g]$ defined on $M$ is a \emph{SMMS invariant} if for all diffeomorphisms $\psi\colon M\to M$,
\[ \psi^\ast\left(T[g,v^m\dvol_g]\right) = T[\psi^\ast g,(\psi^\ast v)^m\dvol_{\psi^\ast g}] . \]
\end{defn}

\begin{remark}

It is through this notion that Chang, Gursky and Yang~\cite{CGY0} related conformal invariants of a SMMS to diffeomorphism invariants.  More precisely, Moser~\cite{Moser1965} showed that given any smooth measure $d\omega$ on a compact Riemannian manifold $(M^n,g)$ such that $\int d\omega=\Vol(M)$, there exists a diffeomorphism such that $\dvol_g=\psi^\ast d\omega$.
\end{remark}

When we wish to use ideas from conformal geometry to study SMMS, it will also be useful to have a definition of a conformal invariant on a SMMS.

\begin{defn}
Let $(M^n,g,v^m)$ be a SMMS.  A SMMS invariant $T[g,v^m\dvol_g]$ on $M$ is said to be \emph{conformally invariant of weight $w$} if for all positive $u\in C^\infty(M)$,
\[ T\left[u^{-2}g,u^{-m-n}v^m\dvol_g\right] = u^{-(m+n)w}T\left[g,v^m\dvol_g\right] . \]
\end{defn}

\begin{remark}
This definition makes sense in the limit $\lv m\rv=\infty$, where our convention is such that a conformal invariant of weight $w$ satisfies
\[ T\left[g,e^{-f-\phi}\dvol_g\right] = e^{-wf}T\left[g,e^{-\phi}\dvol_g\right] \]
for all $f\in C^\infty(M)$.  Note, however, that this convention differs from the usual definition in conformal geometry by a factor of $m+n$.
\end{remark}

Finally, we wish to point out that many problems involving measures on Riemannian manifolds ask that the measure be allowed to degenerate on the boundary; that is, it is common to study a SMMS $(M^n,g,v^m\dvol)$ for $M^n$ a compact manifold with boundary and $v\in C^\infty(M)$ a nonnegative function such that $\partial M=v^{-1}(0)$.  While we specifically rule this out in our definition of a SMMS, many of the ideas we introduce here can nevertheless be (partially) made sense of when boundaries of this type are allowed (cf.\ Section~\ref{sec:defn_qe_simple}; \cite{HePetersenWylie2010}).               
\section{Quasi-Einstein Metrics}
\label{sec:quasi_einstein}

Let us now discuss quasi-Einstein SMMS.  In order to so, we must first describe the appropriate weighted analogues of the Ricci and scalar curvatures on a SMMS.  The first definition is standard when $m>0$, being given by the ($m$-)Bakry-\'Emery Ricci tensor.  The second definition agrees with Perelman's scalar curvature in the case $\lv m\rv=\infty$, and is introduced so that quasi-Einstein metrics admit the natural variational characterization in terms of the total weighted scalar curvature functional.  Finally, we discuss some basic but important results for quasi-Einstein SMMS, especially as they will be needed in the sequel~\cite{Case2010b}.

\subsection{The Curvature of a SMMS}
\label{sec:curvature}

\begin{defn}
Let $(M^n,g,v^m\dvol)$ be a SMMS.  The \emph{Bakry-\'Emery Ricci tensor $\Ric_\phi^m$} is the symmetric $(0,2)$-tensor
\[ \Ric_\phi^m := \Ric - mv^{-1}\nabla^2 v = \Ric + \nabla^2\phi - \frac{1}{m}d\phi\otimes d\phi . \]
The \emph{weighted scalar curvature $R_\phi^m$} is the function
\[ R_\phi^m := R - 2mv^{-1}\Delta v - m(m-1)\lv\nabla v\rv^2 = R + 2\Delta\phi - \frac{m+1}{m}\lv\nabla\phi\rv^2 . \]
\end{defn}

In particular, in keeping with the second goal of this article, we see that the Bakry-\'Emery Ricci tensor is continuous in $m\in\bR\cup\{\pm\infty\}$.

It is important to note that it is generally not the case that the weighted scalar curvature is the trace of the Bakry-\'Emery Ricci tensor.  Indeed, it holds that
\begin{equation}
\label{eqn:scalar_as_trace}
R_\phi^m = \tr\Ric_\phi^m + \Delta_\phi\phi .
\end{equation}
To recover the usual notion of the scalar curvature as the trace of the Ricci curvature, one must instead work in the \auxiliarymanifold.  More precisely, suppose that $(M^n,g,v^m\dvol)$ is a SMMS with $m\in\bN$.  A straightforward computation (cf.\ \cite{Besse,Case_Shu_Wei}) shows that for any $X\in\kX(M)$, it holds that
\[ \Ric_{\og}(X,X) = \Ric_\phi^m(X,X), \]
where $X$ on the left hand side denotes the lift of $X$ to $M\times F$.  On the other hand, one can also easily check that
\[ R_{\og} = R_\phi^m + mv^{-2}R_h . \]
In particular, if one constructs the \auxiliarymanifold\ by taking a warped product with a scalar flat manifold, $R_\phi^m$ is the scalar curvature of the \auxiliarymetric.

While the Bianchi identity doesn't quite hold for the Bakry-\'Emery Ricci tensor, it is not too far off.

\begin{prop}
\label{prop:bianchi_identity}
Let $(M^n,g,v^m\dvol)$ be a SMMS.  Then
\begin{equation}
\label{eqn:bianchi_identity}
\delta_\phi\Ric_\phi^m = \frac{1}{2}dR_\phi^m - \frac{1}{m}\Delta_\phi\phi\,d\phi .
\end{equation}
\end{prop}

\begin{proof}

Straightforward computations show that
\begin{align*}
\delta\Ric_\phi^m & = \frac{1}{2}dR + \Ric(\nabla\phi) + d\Delta\phi - \frac{1}{m}\Delta\phi\,d\phi - \frac{1}{2m}d\lv\nabla\phi\rv^2 \\
\Ric_\phi^m(\nabla\phi) & = \Ric(\nabla\phi) + \frac{1}{2}d\lv\nabla\phi\rv^2 - \frac{1}{m}\lv\nabla\phi\rv^2\,d\phi .
\end{align*}
Taking the difference immediately yields~\eqref{eqn:bianchi_identity}.
\end{proof}

Since in general we know that the weighted scalar curvature is not the trace of the Bakry-\'Emery Ricci tensor, it is also useful to recast Proposition~\ref{prop:bianchi_identity} as the failure of $\Ric_\phi^m$ to be in the kernel of the weighted Bianchi operator $\mB_\phi$, defined by
\[ \mB_\phi T = \delta_\phi T - \frac{1}{2}d\tr T \]
for all $T\in S^2T^\ast M$.

\begin{prop}
\label{prop:bianchi_operator}
Let $(M^n,g,v^m\dvol)$ be a SMMS.  Then
\begin{equation}
\label{eqn:bianchi_ric}
\mB_\phi\Ric_\phi^m = \frac{1}{2}e^{\frac{2}{m}\phi} d\left(e^{-\frac{2}{m}\phi}\Delta_\phi\phi\right) .
\end{equation}
\end{prop}

\begin{proof}

This follows immediately from~\eqref{eqn:scalar_as_trace} and~\eqref{eqn:bianchi_identity}.
\end{proof}

Finally, it will be useful to know how the Bakry-\'Emery Ricci tensor and the weighted scalar curvature transform for conformal changes of a SMMS.

\begin{prop}
\label{prop:curvature_conformal}
Let $(M^n,g,v^m\dvol)$ be a SMMS and let $u\in C^\infty(M)$ be a positive function.  The Bakry-\'Emery Ricci curvature $\Ric_{f,\phi}^m$ and the weighted scalar curvature $u^2R_{f,\phi}^m$ of the SMMS $(M^n,\hat g,\hat v^m\dvol_{\hat g})$ determined by $u$ via~\eqref{eqn:scms} are given by
\begin{align*}
\Ric_{f,\phi}^m & = \Ric_\phi^m + (m+n-2)u^{-1}\nabla^2 u + \left(u^{-1}\Delta_\phi u - (m+n-1)u^{-2}\lv\nabla u\rv^2\right)g \\
R_{f,\phi}^m & = R_\phi^m + 2(m+n-1)u^{-1}\Delta_\phi u - (m+n)(m+n-1)u^{-2}\lv\nabla u\rv^2 .
\end{align*}
\end{prop}

\begin{remark}

The use of the symbol $f$ in denoting $\Ric_{f,\phi}^m$ and $R_{f,\phi}^m$ is to emphasize that these tensors are well-defined when $\lv m\rv=\infty$, as is easily checked.
\end{remark}

\begin{proof}

By the definition of the Bakry-\'Emery Ricci tensor and the weighted scalar curvature, we see that
\begin{align*}
\Ric_{f,\phi}^m & = \Ric_{u^{-2}g} - m\left(vu^{-1}\right)^{-1}\nabla_{u^{-2}g}^2(vu^{-1}) \\
u^2R_{f,\phi}^m & = R_{u^{-2}g} - 2m(vu^{-1})\Delta_{u^{-2}g}(vu^{-1}) - m(m-1)(vu^{-1})^{-2}\lv\nabla (vu)^{-1}\rv_{u^{-2}g}^2 .
\end{align*}
On the other hand, it is well-known (see~\cite{Besse}) that for any $w\in C^\infty(M)$,
\begin{align*}
\Ric_{u^{-2}g} & = \Ric_g + (n-2)u^{-1}\nabla_g^2 u + \left(u^{-1}\Delta_g u - (n-1)u^{-2}\lv\nabla u\rv^2\right)g \\
\nabla_{u^{-2}g}^2 w & = \nabla_g^2 w + u^{-1}dw\otimes du + u^{-1}du\otimes dw - u^{-1}\lp\nabla u,\nabla w\rp_g\,g .
\end{align*}
The result then follows via a straightforward computation.
\end{proof}
\subsection{Quasi-Einstein SMMS}
\label{sec:defn_qe_simple}

\begin{defn}
A SMMS $(M^n,g,v^m\dvol_g)$ is said to be \emph{quasi-Einstein} if there is a constant $\lambda\in\bR$ such that
\begin{equation}
\label{eqn:qe_ric}
\Ric_\phi^m = \lambda g .
\end{equation}
When this is the case, we call $\lambda$ the \emph{quasi-Einstein constant}.

If $(M^n,g,v^m\dvol_g)$ is quasi-Einstein with quasi-Einstein constant $\lambda=0$, we say that it is \emph{BER-flat}.
\end{defn}

An important fact about quasi-Einstein SMMS, established by D.-S.\ Kim and Y.\ H.\ Kim~\cite{Kim_Kim} in the case $0<m<\infty$ and by Hamilton~\cite{Hamilton_1995} in the case $\lv m\rv=\infty$, is that they determine another important constant as a consequence of the Bianchi identity.

\begin{lem}
\label{lem:bianchi}
Let $(M^n,g,e^{-\phi}\dvol,m)$ be a quasi-Einstein SMMS with quasi-Einstein constant $\lambda$.  Then there are constants $\mu,\mu^\prime\in\bR$ such that
\begin{subequations}
\label{eqn:bianchi}
\begin{equation}
\label{eqn:bianchi_finite}
R_\phi^m + m\mu v^{-2} = (m+n)\lambda
\end{equation}
when $\lv m\rv<\infty$ and
\begin{equation}
\label{eqn:bianchi_infinite}
R_\phi^\infty + 2\lambda(\phi-n) = -\mu^\prime
\end{equation}
when $m=\infty$.
\end{subequations}
\end{lem}

\begin{proof}

A direct computation shows that
\[ \mB_\phi(\lambda g) = -\lambda\,d\phi = \frac{m}{2}e^{\frac{2}{m}\phi}d\left(\lambda e^{-\frac{2}{m}\phi}\right) . \]
It then follows from~\eqref{eqn:bianchi_ric}, the assumption $\Ric_\phi^m=\lambda g$, and the connectedness of $M$ that there is a constant $\mu$ such that
\[ \Delta_\phi\phi = m\lambda + m\mu e^{\frac{2}{m}\phi} . \]
The result then follows from the fact
\[ R_\phi^m = \tr\Ric_\phi^m + \Delta_\phi\phi . \qedhere \]
\end{proof}

When $m\in\bN$, the constant $\mu$ has an important implication for constructing auxiliary manifolds (cf.\ \cite{Besse,Case_Shu_Wei,Kim_Kim}).

\begin{cor}
\label{cor:auxiliarymanifold_einstein}
Let $(M^n,g,v^m\dvol)$ be a quasi-Einstein constant with $m\in\bN$ and quasi-Einstein constant $\lambda$, and let $\mu$ be the constant such that~\eqref{eqn:bianchi_finite} holds.  Then given any Riemannian manifold $(F^m,h)$ such that $\Ric_h=\mu h$, the \auxiliarymanifold
\[ \left( M^n\times F^m, \og:=g\oplus v^2h \right) \]
is an Einstein manifold with $\Ric_{\og}=\lambda\og$.
\end{cor}

\begin{remark}
\label{rk:m=1_trouble}
Note that there do exist quasi-Einstein SMMS with $m=1$ and \charconstant\ $\mu\not=0$, despite the fact that there are no one-dimensional Einstein manifolds with nonzero scalar curvature.  It is for this reason that we have phrased Corollary~\ref{cor:auxiliarymanifold_einstein} in this way.

Note also that, since they are mostly interested in quasi-Einstein SMMS as a way to study warped product Einstein manifolds, some authors impose the additional constraint that any quasi-Einstein SMMS with $m=1$ have \charconstant\ $\mu=0$; see, for example, \cite{Corvino2000}.
\end{remark}

It is clear that~\eqref{eqn:bianchi} is continuous as an expression for $m\in[0,\infty)$.  To understand the limiting behavior as $\lv m\rv\to\infty$, the following elementary calculus lemma will be useful.

\begin{lem}
\label{lem:calculus_lemma}
Let $M$ be a smooth manifold and let $f\in C^\infty(M)$.  For any $p\in M$,
\[ \lim_{m\to\infty} m\left(1-\exp\left(-\frac{f(p)}{m}\right)\right) = f(p) . \]
\end{lem}

\begin{prop}
\label{prop:bianchi_limit}
Let $(M^n,g_i,e^{-\phi_i}\dvol_{g_i},m_i)$ be a smooth family of quasi-Einstein SMMS with quasi-Einstein constant $\lambda_i$ and $\lv m_i\rv<\infty$ and let $\mu_i$ be as in~\eqref{eqn:bianchi_finite}.  Suppose additionally that $(M^n,g_i,e^{-\phi_i}\dvol_{g_i},m_i)$ converges to the quasi-Einstein SMMS $(M^n,g,e^{-\phi}\dvol_g,\infty)$ with quasi-Einstein constant $\lambda$, and let $\mu^\prime$ be as in~\eqref{eqn:bianchi_infinite}.  Then
\begin{align*}
\lambda & = \lim_{i\to\infty} \lambda_i = \lim_{i\to\infty} \mu_i \\
\mu^\prime - n\lambda & = \lim_{i\to\infty} m_i\left(\mu_i-\lambda_i\right) .
\end{align*}
\end{prop}

\begin{proof}

By the assumptions on $(M^n,g_i,e^{-\phi_i}\dvol_{g_i},m_i)$ and $(M,g,e^{-\phi}\dvol_g,\infty)$, it suffices to consider the limiting behavior of~\eqref{eqn:bianchi_finite} as $i\to\infty$.  First, dividing by $m_i$, we see that it must be the case that
\[ \lim_{i\to\infty} (\lambda_i-\mu_i) = 0 . \]
On the other hand, the quasi-Einstein equation~\eqref{eqn:qe_ric} implies that $\lambda_i\to\lambda$ as $i\to\infty$, yielding the first claim.  Next, we can rewrite~\eqref{eqn:bianchi_finite} as
\[ R_\phi^m + m\mu_i\left(e^{\frac{2}{m_i}\phi_i} - 1\right) - 2n\lambda_i = m_i(\lambda_i-\mu_i) - n\lambda_i . \]
The result then follows from Lemma~\ref{lem:calculus_lemma} by taking the limit as $i\to\infty$.
\end{proof}

The constant $\mu$ plays an important role in the study of quasi-Einstein SMMS, and for this reason it will be convenient to give it a name.

\begin{defn}
Let $(M^n,g,v^m\dvol)$ be a quasi-Einstein SMMS with quasi-Einstein constant $\lambda$.  The \emph{\charconstant} $\mu$ is the constant such that~\eqref{eqn:bianchi} holds.
\end{defn}

By Proposition~\ref{prop:bianchi_limit}, this convention is such that the \charconstant\ and the quasi-Einstein constant coincide for quasi-Einstein SMMS with $\lv m\rv=\infty$.

It will be important for us to understand what it means for a SMMS to be conformally quasi-Einstein.  In this way, the following consequence of Proposition~\ref{prop:curvature_conformal} will be useful.

\begin{prop}
\label{prop:curvature_conformal_qe}
Let $(M^n,g,v^m\dvol)$ be a quasi-Einstein SMMS and suppose that $u\in C^\infty(M)$ is such that the SMMS $(M^n,\hat g,\hat v^m\dvol_{\hat g})$ determined by~\eqref{eqn:scms} is quasi-Einstein with quasi-Einstein constant $\lambda$ and \charconstant\ $\mu$.  Then
\begin{subequations}
\label{eqn:qe_uv}
\begin{align}
\label{eqn:qe_tf_uv} 0 & = \left(uv\Ric + (m+n-2)v\nabla^2u - mu\nabla^2v\right)_0 \\
\label{eqn:qe_lambda_uv} n\lambda v^2 & = (uv)^2R + (m+2n-2)uv^2\Delta u - mu^2v\Delta v \\
\notag & \quad - (m+n-1)nv^2\lv\nabla u\rv^2 + mnuv\lp\nabla u,\nabla v\rp \\
\label{eqn:qe_mu_uv} n\mu u^2 & = (uv)^2R + (m+n-2)uv^2\Delta u - (m-n)u^2v\Delta v \\
\notag & \quad - (m+n-2)nuv\lp\nabla u,\nabla v\rp + (m-1)nu^2\lv\nabla v\rv^2 ,
\end{align}
\end{subequations}
where $T_0=T-\frac{1}{n}\tr_g T\,g$ denotes the tracefree part of a section $T\in S^2T^\ast M$ and all inner products and traces on the right hand side are computed with respect to $g$.
\end{prop}

\begin{proof}

Let $\Ric_{f,\phi}^m$ and $u^2R_{f,\phi}^m$ denote the Bakry-\'Emery Ricci tensor and the weighted scalar curvature of $(M^n,\hat g,\hat v^m\dvol_{\hat g})$ as in Proposition~\ref{prop:curvature_conformal}.  The tracefree part of the condition $0=\Ric_{f,\phi}^m-\lambda\hat g$ is precisely~\eqref{eqn:qe_tf_uv}, while its trace is~\eqref{eqn:qe_lambda_uv}.  By definition, $\mu$ is determined by the equation
\[ u^2R_{f,\phi}^m + m\mu(vu^{-1})^{-2} = (m+n)\lambda . \]
\eqref{eqn:qe_mu_uv} then follows by using~\eqref{eqn:qe_lambda_uv} to rewrite $\lambda$ in terms of $u$ and $v$.
\end{proof}

There are three important observations to make about Proposition~\ref{prop:curvature_conformal_qe}.

First, \eqref{eqn:qe_uv} make sense when $u$ and $v$ are allowed to be zero or change signs.  The zero sets of such functions play an important role in various geometric problems.  For example, when $m=0$ the zero set of $u$ corresponds to the conformal infinity of a Poincar\'e-Einstein metric (cf.\ \cite{Gover2008}), while when $m=1$ the zero set of $v$ corresponds to a horizon of a static metric (cf.\ \cite{Corvino2000,HePetersenWylie2010}).

Second, \eqref{eqn:qe_lambda_uv} and~\eqref{eqn:qe_mu_uv} suggest that the constants $\lambda$ and $\mu$ be thought of as the squared-lengths of $u$ and $v$, respectively.  More precisely, given constants $c,k>0$, one readily verifies that~\eqref{eqn:qe_tf_uv} is invariant under the rescalings $u\mapsto cu$ and $v\mapsto kv$, while~\eqref{eqn:qe_lambda_uv} and~\eqref{eqn:qe_mu_uv} imply that $\lambda\mapsto c^2\lambda$ and $\mu\mapsto k^2\mu$.  For this reason, it will sometimes be useful to consider a \emph{SMMS with \charconstant\ $\mu$}, which is a SMMS together with a fixed choice of $\mu$ to be used in defining natural geometric objects on the SMMS.  For example, we make the following definition of a quasi-Einstein scale.

\begin{defn}
Let $(M^n,g,v^m\dvol)$ be a SMMS with \charconstant\ $\mu$.  A \emph{quasi-Einstein scale} is a function $u\in C^\infty(M)$ such that the SMMS $(M^n,\hat g,\hat v^m\dvol_{\hat g})$ determined by~\eqref{eqn:scms}, wherever it is defined, is a quasi-Einstein SMMS with \charconstant\ $\mu$.

If $u$ is a quasi-Einstein scale, the \emph{quasi-Einstein constant (of $u$)} is the quasi-Einstein constant of $(M^n,\hat g,\hat v^m\dvol_{\hat g})$.
\end{defn}

Note that in this definition we do not impose the constraint $u>0$.

Third, Proposition~\ref{prop:curvature_conformal_qe} reveals a symmetry present in the study of conformally quasi-Einstein SMMS which gives rise to the useful notion of duality.

\begin{cor}
\label{cor:duality}
Let $(M^n,g,v^m\dvol)$ be a SMMS with \charconstant\ $\mu$ and suppose that $u\in C^\infty(M)$ is a positive quasi-Einstein scale with quasi-Einstein constant $\lambda$.  Then the SMMS $(M^n,g,u^{2-m-n}\dvol)$ with \charconstant\ $\lambda$ is such that $v\in C^\infty(M)$ is a quasi-Einstein scale with quasi-Einstein constant $\mu$.
\end{cor}

\begin{proof}

This follows immediately from the observation that~\eqref{eqn:qe_uv} is invariant under the change of variables
\[ (u,v,\lambda,\mu,m) \mapsto (v,u,\mu,\lambda,2-m-n) . \qedhere \]
\end{proof}

Indeed, the duality observed in Corollary~\ref{cor:duality} gives exactly the duality between two methods discussed in Section~\ref{sec:warped} and Section~\ref{sec:conformal} for studying smooth metric measure spaces: The approach of Bakry-\'Emery is to study SMMS $(M^n,g,v^m\dvol)$ with $m\geq 0$ as we have defined here, while the approach of Chang-Gursky-Yang is to study the SMMS $(M^n,g,u^{2-m-n}\dvol_g)$ by studying the scale (i.e.\ conformal factor) $u\in C^\infty(M)$ on the SMMS $(M^n,g,1^m\dvol_g)$, again with $m\geq 0$.

An equivalent formulation of Corollary~\ref{cor:duality} which is perhaps more useful in practice is the following.

\begin{cor}
\label{cor:duality_no_scales}
Let $(M^n,g)$ be a Riemannian manifold and fix $m\in\bR\cup\{\pm\infty\}$ and positive functions $u,v\in C^\infty(M)$.  Then the following are equivalent:
\begin{enumerate}
\item $(M^n,g,v^m\dvol_g)$ is a SMMS with \charconstant\ $\mu$ such that $u$ is a quasi-Einstein scale with quasi-Einstein constant $\lambda$.
\item $(M^n,\hat g,\hat v^m\dvol_{\hat g})$ is a quasi-Einstein SMMS with quasi-Einstein constant $\lambda$ and \charconstant\ $\mu$, where $\hat g=u^{-2}g$ and $\hat v=v/u$.
\item $(M^n,g,u^{2-m-n}\dvol_g)$ is a SMMS with \charconstant\ $\lambda$ such that $v$ is a quasi-Einstein scale with quasi-Einstein constant $\mu$.
\item $(M^n,\tilde g,\tilde u^{2-m-n}\dvol_{\tilde g})$ is a quasi-Einstein SMMS with quasi-Einstein constant $\mu$ and \charconstant\ $\lambda$, where $\tilde g=v^{-2}g$ and $\tilde u=u/v$.
\end{enumerate}
\end{cor}

\begin{proof}

The equivalence of the first and second conditions, and likewise of the third and fourth conditions, is due to the definition of a quasi-Einstein scale.  The equivalence of the first and third conditions is due to Corollary~\ref{cor:duality}.
\end{proof}
\subsection{Variational Properties}
\label{sec:variational}

It is well-known that Einstein metrics arise as critical points of the total scalar curvature functional.  Similarly, Perelman~\cite{Perelman1} showed that shrinking and steady gradient Ricci solitons arise as critical points of his $\mW$ and $\mF$-functionals, respectively, and later Feldman, Ilmanen and Ni~\cite{FeldmanIlmanenNi2005} showed that expanding gradient Ricci solitons arise as critical points of their $\mW^+$-functional.  As generalizations of Einstein metrics and gradient Ricci solitons, it is natural to expect a similar observation to hold for quasi-Einstein metrics.  In this section, we verify this fact, which also illustrates the relationship between the total scalar curvature functional and the $\mF$, $\mW$, and $\mW^+$-functionals.

In the following, given a smooth manifold $M^n$, we denote by $\Met(M)$ the space of Riemannian metrics on $M$ and by $\kM$ the space of smooth measures on $M$.

\begin{defn}
Let $M^n$ be a smooth manifold and fix constants $m\in\bR\cup\{\pm\infty\}$ and $\mu\in\bR$.  The \emph{$(m,\mu)$-energy functional $\mW_\mu^m\colon\Met(M)\times\kM\to\bR$} is given by
\begin{subequations}
\label{eqn:mW}
\begin{equation}
\label{eqn:mW_finite}
\mW_\mu^m\left(g,v^m\dvol_g\right) = \int_M \left( R_\phi^m + m\mu v^{-2}\right) v^m\dvol_g
\end{equation}
when $\lv m\rv<\infty$ and by
\begin{equation}
\label{eqn:mW_infinite}
\mW_\mu^m\left(g,e^{-\phi}\dvol_g\right) = \int_M \left( R_\phi^m + 2\mu(\phi-n)\right) e^{-\phi}\dvol_g
\end{equation}
\end{subequations}
when $\lv m\rv=\infty$.

When the context is clear, we shall often refer to the $(m,\mu)$-energy functional as the \emph{energy functional}.
\end{defn}

In other words, the $(m,\mu)$-energy functional assigns to each pair $(g,v^m\dvol)\in\Met(M)\times\kM$ the \emph{total weighted scalar curvature of the SMMS $(M^n,g,v^m\dvol)$ with \charconstant\ $\mu$}.  As special cases, observe that when $m=0$, the energy functional is the usual total scalar curvature functional in Riemannian geometry, while when $\lv m\rv=\infty$, it is, up to a multiplicative constant, equal to one of the entropy functionals appearing in the study of the Ricci flow~\cite{FeldmanIlmanenNi2005,Perelman1}, depending upon the sign of $\mu$.  More precisely, recall that Perelman introduced the functionals $\mF$ and $\mW$~\cite{Perelman1} and Feldman, Ilmanen and Ni~\cite{FeldmanIlmanenNi2005} introduced the functional $\mW^+$ as
\begin{align*}
\mF(g,\phi) & = \int_M R_\phi^\infty\,e^{-\phi}\dvol_g \\
\mW(g,\phi,\tau) & = \int_M \tau\left( R_\phi^\infty + \tau^{-1}(\phi-n)\right) (4\pi\tau)^{-\frac{n}{2}} e^{-\phi}\dvol_g \\
\mW^+(g,\phi,\tau) & = \int_M \tau\left( R_\phi^\infty - \tau^{-1}(\phi-n)\right) (4\pi\tau)^{-\frac{n}{2}} e^{-\phi}\dvol_g
\end{align*}
for any $g\in\Met(M)$, $f\in C^\infty(M)$ and $\tau>0$.  In this way, direct comparison with~\eqref{eqn:mW_infinite} yields that the $(\infty,\mu)$-energy functional corresponds to the $\mW^+$, $\mF$, and $\mW$-functionals when $\mu$ is negative, zero, or positive, respectively, and with $\lv\mu\rv$ equivalent to the scale $\tau$.

Second, modulo a suitable renormalization, the energy functional depends continuously on its parameters, including the dimensional parameter $m\in\bR\cup\{\pm\infty\}$.  More precisely, we have the following simple proposition.

\begin{prop}
\label{prop:energy_functional_smooth}
Let $(M^n,g)$ be a compact Riemannian manifold and fix $\phi\in C^\infty(M)$, $\mu\in\bR$.  Then
\[ \lim_{m\to\infty} \mW_\mu^m\left(g,e^{-\phi}\dvol,m\right) - (m+2n)\Vol_\phi(M) = \mW_\mu^\infty\left(g,e^{-\phi}\dvol,\infty\right) . \]
\end{prop}

\begin{proof}

This follows as in the proof of Proposition~\ref{prop:bianchi_limit}
\end{proof}

In particular, Proposition~\ref{prop:energy_functional_smooth} can be regarded as stating that the energy functional interpolates between the total scalar curvature functional and the entropy functionals appearing in the Ricci flow.

To see that the $(m,\mu)$-energy functionals yield a variational characterization of quasi-Einstein metrics, we first need to compute the first variation of $\mW_\mu^m$.

\begin{prop}
\label{prop:gen_var}
Let $M^n$ be a smooth manifold and fix constants $m\in\bR\cup\{\pm\infty\}$ and $\mu\in\bR$.  Let $(g,e^{-\phi}\dvol_g)\in\Met(M)\times\kM$ and let $(h,\psi)=(\delta g,\delta\phi)$ be a compactly supported variation of $g$ and $h$.  Then
\begin{align*}
\delta\mW_\mu^m & = -\int_M \bigg[ \lp\Ric_\phi^m - \frac{1}{2}(R_\phi^m+m\mu v^{-2})g,h\rp \\
& \qquad + \left(R_\phi^m - \frac{2}{m}\Delta_\phi\phi + (m-2)\mu v^{-2}\right)\psi \bigg] e^{-\phi}\dvol_g ,
\end{align*}
where all derivatives and inner products are taken with respect to $g$ and we define
\[ \delta\mW_\mu^m = \frac{d}{ds}\mW_\mu^m\left(g+sh,e^{-\phi+s\psi}\dvol_{g+sh}\right)\big|_{s=0} . \]
\end{prop}

\begin{proof}

We begin by computing the derivative
\[ \delta R_\phi^m=\frac{d}{ds}R_\phi^m[g+sh,e^{-\phi+s\psi}\dvol]\big|_{s=0} . \]
It is well-known that
\begin{align*}
\delta R & = -\lp\Ric,h\rp + \delta^2 h - \Delta\tr h \\
\delta\Delta\phi & = \Delta\psi - \lp h,\nabla^2\phi\rp - \lp\delta h,d\phi\rp + \frac{1}{2}\lp d\tr h,d\phi\rp \\
\delta\lv\nabla\phi\rv^2 & = 2\lp\nabla\phi,\nabla\psi\rp - \lp h,d\phi\otimes d\phi\rp ,
\end{align*}
where $\delta$ on the right hand side denotes the divergence (see, for example, \cite{Topping2006}).  On the other hand, we easily check that
\begin{align*}
\delta_\phi^2 h & = \delta_\phi\left(\delta h - h(\nabla\phi,\cdot)\right) \\
& = \delta^2 h - 2\lp\delta h,d\phi\rp - \lp h,\nabla^2\phi\rp + \lp h,d\phi\otimes d\phi\rp .
\end{align*}
Combining these, it follows easily that
\begin{equation}
\label{eqn:deltaRphim}
\delta R_\phi^m = -\lp\Ric_\phi^m, h\rp + \delta_\phi^2 h - \Delta_\phi\tr h + 2\left(\Delta_\phi\psi - \frac{1}{m}\lp\nabla\phi,\nabla\psi\rp\right) .
\end{equation}
The result then follows from the formula
\[ \delta\mW_\mu^m = \int_M \left[\delta R_\phi^m + 2\mu v^{-2}\psi - \left(R_\phi^m+m\mu v^{-2}\right)\left(\psi-\frac{1}{2}\tr h\right)\right] e^{-\phi}\dvol \]
and integration by parts.
\end{proof}

\begin{remark}
\label{rk:bianchi_remark}
Proposition~\ref{prop:gen_var} yields an alternative proof of the Bianchi identity~\eqref{eqn:bianchi_identity} via the diffeomorphism invariance of the energy functional (cf.\ \cite{Kazdan1981}).  Explicitly, let $\phi_t$ be a one parameter family of diffeomorphisms generated by a vector field $X$ and consider the corresponding variations of the $(m,0)$-energy functional.  By definition, $\delta g=L_Xg$ and $\delta\phi=X\phi$, and so the diffeomorphism invariance of the energy functional together with Proposition~\ref{prop:gen_var} imply that
\begin{align*}
0 & = - \int_M\left[\lp\Ric_\phi^m-\frac{1}{2}R_\phi^mg,L_Xg\rp + \left(R_\phi^m-\frac{2}{m}\Delta_\phi\phi\right)X\phi\right] e^{-\phi}\dvol \\
& = 2\int_M\lp\delta_\phi\Ric_\phi^m - \frac{1}{2}dR_\phi^m + \frac{1}{m}\Delta_\phi\phi\,d\phi, X\rp e^{-\phi}\dvol .
\end{align*}
\end{remark}

To arrive at the desired variational characterization of quasi-Einstein metrics, we must restrict to variations which lie in
\[ \mC_1(M,m) = \left\{ (g,v^m\dvol_g)\in\Met(M)\times\kM \colon \int_M v^m\dvol_g = 1 \right\} , \]
the space of smooth metrics and measures on a smooth manifold $M^n$ such that the SMMS $(M^n,g,v^m\dvol)$ has unit weighted volume.

\begin{prop}
\label{prop:qe_var}
Let $M^n$ be a compact manifold and fix constants $m\in\bR\cup\{\pm\infty\}$ and $\mu\in\bR$.  Let $(g,e^{-\phi}\dvol_g)\in\mC_1(M,m)$ and consider a variation $(h,\psi)=(\delta g,\delta\phi)$ of $(g,e^{-\phi}\dvol_g)$ which remains in $\mC_1(m)$.  Then
\begin{align}
\label{eqn:constraint_implication} 0 & = \int \left(\psi-\frac{1}{2}\tr_g h\right)e^{-\phi}\dvol_g \\
\label{eqn:wmumvar} \delta\mW_\mu^m & = -\int\bigg[\lp\Ric_\phi^m - \frac{1}{m}\Delta_\phi\phi\,g - \mu v^{-2}g,h\rp \\
\notag & \qquad + \left(R_\phi^m-\frac{2}{m}\Delta_\phi\phi + (m-2)\mu v^{-2}\right)\left(\psi-\frac{1}{2}\tr h\right)\bigg] e^{-\phi}\dvol_g .
\end{align}

In particular, if $(g,v^m\dvol_g)\in\mC_1(M,m)$ is a critical point of the energy functional, then $(M^n,g,v^m\dvol_g)$ is a quasi-Einstein SMMS with \charconstant\ $\mu$.
\end{prop}

\begin{proof}

\eqref{eqn:constraint_implication} follows immediately from the fact
\[ 0 = \frac{d}{ds}\int_M e^{-\phi_s}\dvol_{g_s} , \]
where $(g_s,e^{-\phi_s}\dvol_{g_s})$ is a curve in $\mC_1(M,m)$ with derivative at $s=0$ given by $(h,\psi)$.  \eqref{eqn:wmumvar} follows immediately from Proposition~\ref{prop:gen_var} by rewriting the second summand in the integral as a product with $\psi-\frac{1}{2}\tr h$.  Finally, if $(g,v^m\dvol_g)$ is a critical point of the energy functional, we must have that
\begin{align*}
0 & = \Ric_\phi^m - \frac{1}{m}\Delta_\phi\phi\,g - \mu v^{-2}g \\
(m+n-2)\lambda & = R_\phi^m - \frac{2}{m}\Delta_\phi\phi + (m-2)\mu v^{-2}
\end{align*}
for some constant $\lambda$.  The result then follows by taking the trace of the first equation and using~\eqref{eqn:scalar_as_trace} to compare it to the second equation.
\end{proof}
\subsection{Compact Quasi-Einstein SMMS}
\label{sec:cpt_results}

When studying quasi-Einstein SMMS, a useful ``toy model'' are the compact SMMS, where it is easier to use integral estimates and the maximum principle to establish \emph{a priori} estimates on the quasi-Einstein potential $\phi$ (cf.\ \cite{Case_Shu_Wei,Rimoldi2010}).  Here we present a sample of results of this type, focusing on those which will be useful in the sequel~\cite{Case2010b}.

First, the result of D.-S.\ Kim and Y.\ H.\ Kim~\cite{Kim_Kim} together with the well-known fact that the only nontrivial compact gradient Ricci solitons are shrinkers yields the following result for quasi-Einstein SMMS.

\begin{prop}
\label{prop:compact_positive_constants}
A compact quasi-Einstein SMMS $(M^n,g,e^{-\phi}\dvol_g,m)$ is trivial if the quasi-Einstein constant is nonpositive or if $\lv m\rv<\infty$ and the \charconstant\ is nonpositive.
\end{prop}

The remainder of the results we present rely on results from comparison geometry for SMMS, which in particular requires one to take $m\geq 0$.  In some cases, ``dual'' results can be obtained via Corollary~\ref{cor:duality}, as we shall illustrate in a few examples.

To begin, we have the following generalization of Myers' theorem for SMMS due to Qian~\cite{Qian1997}.

\begin{prop}
\label{prop:myers}
Let $(M^n,g,v^m\dvol)$ be a quasi-Einstein SMMS with $m<\infty$ and positive quasi-Einstein constant.  Then $M$ is compact.
\end{prop}

Note that this result is false in the limit $m\to\infty$, as evidenced by the Gaussian shrinker; see Section~\ref{sec:examples} for details.

Next, recall that for shrinking gradient Ricci solitons, one has a gradient estimate for the potential which establishes that the potential grows at most quadratically with respect to the metric distance (cf.\ \cite{Cao2009_surveyb}).  This estimate can be generalized to quasi-Einstein SMMS with $m>1$ or $m<1-n$, and is closely related to the following scalar curvature estimate proven in~\cite{Case_Shu_Wei}.

\begin{prop}
\label{prop:scalar_bd_warped}
Let $(M^n,g,v^m\dvol_g)$ be a compact quasi-Einstein SMMS with $m>1$ and positive quasi-Einstein constant $\lambda$.  Then the scalar curvature $R$ satisfies
\[ R > \frac{n(n-1)}{m+n-1}\lambda . \]
\end{prop}

\begin{remark}

The proof given in~\cite{Case_Shu_Wei} in fact allows one to conclude that if $(M^n,g,v^m\dvol_g)$ is a compact quasi-Einstein SMMS with boundary $v^{-1}(0)$, $m>1$, and positive quasi-Einstein constant, then
\[ R \geq \frac{n(n-1)}{m+n-1} . \]
\end{remark}

As an immediate corollary, we have the following \emph{a priori} gradient estimate for the quasi-Einstein potential of a nontrivial compact quasi-Einstein SMMS.

\begin{cor}
\label{cor:gradient_est_bigger1}
Let $(M^n,g,v^m\dvol_g)$ be a compact quasi-Einstein SMMS with $m>1$, quasi-Einstein constant $\lambda>0$ and \charconstant\ $\mu>0$.  Then
\begin{equation}
\label{eqn:dil}
\lv\nabla v\rv^2 + \frac{\lambda}{m+n-1}v^2 < \frac{\mu}{m-1} .
\end{equation}
\end{cor}

\begin{proof}

Taking the trace of the quasi-Einstein equation~\eqref{eqn:qe_ric} and using~\eqref{eqn:bianchi}, we see that
\[ (n-m)\lambda = R + m(m-1)v^{-2}\lv\nabla v\rv^2 - m\mu v^{-2} . \]
The result then follows from Proposition~\ref{prop:scalar_bd_warped}.
\end{proof}

By duality, a similar gradient estimate holds for quasi-Einstein SMMS with dimensional parameter $m<1-n$.

\begin{cor}
\label{cor:dil}
Let $(M^n,g,u^{2-m-n}\dvol_g)$ be a compact quasi-Einstein SMMS with $m>1$, quasi-Einstein constant $\mu>0$, and \charconstant\ $\lambda$.  Then
\begin{equation}
\label{eqn:dil_estimate}
\lv\nabla u\rv^2 + \frac{\lambda}{m+n-1} < \frac{\mu}{m-1}u^2 .
\end{equation}
\end{cor}

\begin{proof}

By Corollary~\ref{cor:duality_no_scales}, it follows that the SMMS
\[ \left( M^n, u^{-2}g, (u^{-1})^m\dvol_{u^{-2}g} \right) \]
is a quasi-Einstein SMMS with quasi-Einstein constant $\lambda$ and \charconstant\ $\mu$.  Applying Corollary~\ref{cor:gradient_est_bigger1} then yields~\eqref{eqn:dil_estimate}.
\end{proof}

Duality also yields the following analogue of Proposition~\ref{prop:scalar_bd_warped} for quasi-Einstein SMMS with dimensional parameter $m<1-n$.

\begin{cor}
\label{cor:scal_bd}
Let $(M^n,g,u^{2-m-n}\dvol_g)$ be a compact quasi-Einstein SMMS with $m>1$, quasi-Einstein constant $\mu$, and \charconstant\ $\lambda>0$.  Then the scalar curvature $R_g$ satisfies
\begin{equation}
\label{eqn:scal_bd}
-\frac{n(n-1)}{m-1}\mu < R \leq (m+2n-2)\mu .
\end{equation}
\end{cor}

\begin{proof}

Taking the trace of~\eqref{eqn:qe_ric} and combining it with~\eqref{eqn:bianchi} yields
\begin{equation}
\label{eqn:r_est}
R + (m+n-1)(m+n-2)u^{-2}\lv\nabla u\rv^2 + (m+n-2)\lambda u^{-2} = (m+2n-2)\mu .
\end{equation}
The positivity of $\lambda$ immediately yields the upper bound on $R$, while the lower bound follows immediately from~\eqref{eqn:dil_estimate}.
\end{proof}

\begin{remark}
In~\cite{Case2010b}, we establish a variant on this result which provides an upper bound which is independent of $m$.  More precisely, we show that $\sup\{R,f,\lv\nabla f\rv^2\}\leq ar^2+b$ for constants $a,b>0$ and $r$ the distance from a minimizer of $u=e^{\frac{f}{m+n-2}}$, recovering the well-known result for gradient Ricci solitons.
\end{remark}       
\section{Examples}
\label{sec:examples}

We conclude this article by considering some examples of quasi-Einstein SMMS.  The examples we consider are all well-known in the cases $m\in\bN\cup\{\infty\}$, and our only claim to originality is to illustrate how these examples fit into one-parameter families of quasi-Einstein SMMS parameterized by the dimensional parameter $m$.

The examples we discuss are as follows.  First, we begin our examples with a few comments on trivial quasi-Einstein SMMS.  Second, we discuss the model quasi-Einstein metrics, which we term elliptic Gaussians, and discuss in what sense they are models.  Third, we recast in the language of SMMS the Ricci flat warped products over rotationally symmetric bases~\cite{Berard-Bergery1982,Bohm1999} and their counterparts as gradient Ricci solitons~\cite{Chow_etal,Hamilton1986b}.  Fourth, we say some words about gradient Ricci solitons and Einstein warped products over bases which are $S^2$-bundles over K\"ahler-Einstein manifolds~\cite{Berard-Bergery1982,Cao1996,ChaveValent1996,Koiso1990,LPP,Page1979,PagePope1987}.  Finally, we discuss two simple product constructions which can be used to produce new quasi-Einstein SMMS.

\subsection{Einstein Manifolds}
\label{sec:einstein}

The most basic examples of quasi-Einstein SMMS are the \emph{trivial} examples; i.e.\ SMMS $(M^n,g,1^m\dvol_g)$ for which $(M^n,g)$ is an Einstein manifold.  These examples are quasi-Einstein for any choice of dimensional parameter $m$, but they are also by definition the only quasi-Einstein SMMS with $m=0$.  For this reason, it is natural to think of the notion of a quasi-Einstein SMMS as interpolating between Einstein metrics and other quasi-Einstein metrics.

In order to make comparisons with other quasi-Einstein SMMS, it is useful to discuss the \charconstant\ of a trivial quasi-Einstein SMMS.

\begin{prop}
Let $(M^n,g)$ be an Einstein manifold such that $\Ric=\lambda g$, and let $m\in\bR\cup\{\pm\infty\}$ be a fixed constant.  Then the SMMS $(M^n,g,1^m\dvol_g)$ is a quasi-Einstein SMMS with quasi-Einstein constant $\lambda$ and \charconstant\ $\mu$.
\end{prop}

\begin{proof}

This follows immediately from~\eqref{eqn:bianchi}.
\end{proof}

Of course, we can change our normalization of the measure and consider instead $c^m\dvol_g$ for any constant $c>0$.  While this will change the \charconstant, it will not change its sign.  In particular, the quasi-Einstein constant and the \charconstant\ have the same sign for a trivial quasi-Einstein SMMS.  As we will see, this is not always the case.
\subsection{The Model Spaces: Elliptic Gaussians}
\label{sec:models}

In Section~\ref{sec:cpt_results}, we derived some simple estimates for the scalar curvature and the gradient of the potential of a quasi-Einstein SMMS.  Indeed, these estimates are sharp, and are realized by the \emph{model spaces}, which, following~\cite{HePetersenWylie2010}, we term elliptic Gaussians.  From the standpoint of comparison geometry, special cases of these examples have been observed by Bakry and Qian~\cite{BakryQian2000} in the context of comparison results for the Bakry-\'Emery Ricci tensors.

In order to more clearly discuss these examples as model spaces, it will be useful to state them separately for the cases of positive and negative quasi-Einstein constant.

\begin{defn}
\label{defn:peg}
Fix $m\geq0$ and define $k=\sqrt{m+n-1}$.  The \emph{positive elliptic $m$-Gaussian} is the SMMS
\[ \left(S_+^n, g=dr^2\oplus\left(k\sin(\frac{r}{k})\right)^2d\theta^2, \cos^m(\frac{r}{k})\,\dvol_g\right) ,\]
where $d\theta^2$ denotes the standard metric of constant sectional curvature one on $S^{n-1}$, so that the metric $g$ is the standard metric on $S^n$, normalized so that $\Ric_g=\frac{n-1}{m+n-1}$, and $S_+^n$ is the hemisphere $\{r<\frac{k\pi}{2}\}$.

The \emph{positive elliptic $\infty$-Gaussian} is the SMMS
\[ \left( [0,\infty)\times S^{n-1}, dr^2\oplus r^2d\theta^2, e^{-\frac{r^2}{2}}\dvol_g, \infty\right) . \]
\end{defn}

There are two points to make about our notation.  First, the compactification of the positive elliptic Gaussian is actually a SMMS with boundary with the property that the measure degenerates along the boundary (cf.\ \cite{HePetersenWylie2010}).  Second, our normalization is such that the positive elliptic $m$-Gaussians form a smooth one-parameter family of quasi-Einstein SMMS parameterized by the dimensional parameter $m$.

\begin{prop}
\label{prop:peg_properties}
The positive elliptic $m$-Gaussian is a quasi-Einstein SMMS with quasi-Einstein constant $\lambda=1$ and \charconstant\ $\mu=\frac{m-1}{m+n-1}$.  Moreover, the positive elliptic $m$-Gaussians converge to the positive elliptic $\infty$-Gaussian as $m\to\infty$ in the pointed measured Cheeger-Gromov sense, where we have fixed the base point $r=0$.
\end{prop}

\begin{proof}

It is straightforward to check that
\[ \nabla_g^2\cos\left(\frac{r}{k}\right) = -\frac{1}{k^2}\cos\left(\frac{r}{k}\right)\,g, \]
and so the claim $\Ric_\phi^m=g$ follows from our normalization of $g$.  The computation of the \charconstant\ follows immediately from the definition~\eqref{eqn:bianchi}.  Finally, it is clear from the explicit form of the metric and the measure that the positive elliptic $m$-Gaussians converge to the positive elliptic $\infty$-Gaussian as $m\to\infty$.
\end{proof}

As noted above, one difficulty in working with the positive elliptic Gaussians in general is that the measure degenerates along their boundaries.  This can be remedied by considering their dual spaces; by Corollary~\ref{cor:duality_no_scales}, it follows that the SMMS
\[ \left( H^n, dr^2\oplus\left(k\sinh(\frac{r}{k})\right)^2d\theta^2, 1^m\dvol_g \right) \]
with \charconstant\ $\frac{m-1}{m+n-1}$ admits the function $u=\cosh(\frac{r}{k})$ as a quasi-Einstein scale with quasi-Einstein constant one.  In this form, the positive elliptic Gaussian is a complete SMMS with quasi-Einstein scale $u>0$, both of which are desirable properties.  It is in this sense that we will consider the positive elliptic Gaussian in~\cite{Case2010b}.

We can also see immediately why we refer to the positive elliptic Gaussian as a model: Straightforward computations reveal that, as normalized in Definition~\ref{defn:peg}, the positive elliptic Gaussian has scalar curvature $R=\frac{n(n-1)}{m+n-1}$ and the density $v$ satisfies
\[ \lv\nabla v\rv^2 + \frac{1}{m+n-1}v^2 = \frac{\mu}{m-1} . \]
In other words, equality holds in Proposition~\ref{prop:scalar_bd_warped} and in Corollary~\ref{cor:gradient_est_bigger1} for the elliptic Gaussian.

We also have the corresponding notion of the negative elliptic Gaussian, which is the negatively-curved analogue of the positive elliptic Gaussian.

\begin{defn}
Fix $m\geq0$ and define $k=\sqrt{m+n-1}$.  The \emph{negative elliptic $m$-Gaussian} is the SMMS
\[ \left(H^n, g=dr^2\oplus\left(k\sinh(\frac{r}{k})\right)^2d\theta^2, \cosh^m(\frac{r}{k})\,\dvol_g\right) ,\]
so that $g$ is the standard metric on hyperbolic space, normalized so that $\Ric_g=-\frac{n-1}{m+n-1}$.

The \emph{negative elliptic $\infty$-Gaussian} is the SMMS
\[ \left( \bR^n, dr^2\oplus r^2d\theta^2, e^{\frac{r^2}{2}}\dvol_g, \infty\right) . \]
\end{defn}

Analogous  to Proposition~\ref{prop:peg_properties}, we observe that the negative elliptic Gaussians form a one-parameter family of quasi-Einstein SMMS.

\begin{prop}
\label{prop:neg_properties}
The negative elliptic $m$-Gaussian is a quasi-Einstein SMMS with quasi-Einstein constant $\lambda=-1$ and \charconstant\ $\mu=-\frac{m-1}{m+n-1}$.  Moreover, the negative elliptic $m$-Gaussians converge to the negative elliptic $\infty$-Gaussian as $m\to\infty$ in the pointed measured Cheeger-Gromov sense, where we have fixed the base point $r=0$.
\end{prop}

As formulated, the negative elliptic Gaussian is already a complete quasi-Einstein SMMS, and so we do not need to take its dual to yield an equivalent complete SMMS.  Nevertheless, it is interesting to note that the negative elliptic Gaussian is conformally equivalent to
\[ \left(S_+^n, dr^2\oplus\left(k\sin(\frac{r}{k})\right)^2d\theta^2, 1^m\dvol_g \right) . \]
In other words, the negative elliptic Gaussian is conformally equivalent to a hemisphere in the standard $n$-sphere, normalized so that $\Ric=\frac{n-1}{m+n-1}$, equipped with its standard Riemannian volume element.  In particular, one might study conformally compact SMMS by analogy with the negative elliptic Gaussian.

Finally, we note that both the positive and the negative elliptic Gaussians make sense when $n=1$.  In this form, they serve as the model spaces for the Laplacian comparison theorem~\cite{BakryQian2000}, and form two of the five nontrivial families of one-dimensional quasi-Einstein SMMS (cf.\ \cite{Besse,HePetersenWylie2010}).
\subsection{Rotationally Symmetric Quasi-Einstein SMMS}
\label{sec:cigars}

In this section we discuss examples of complete rotationally symmetric BER-flat SMMS on $\bR^n$ with $m\geq 0$.  These examples have been completely classified through work of B\'erard-Bergery~\cite{Berard-Bergery1982}, Hamilton~\cite{Hamilton1986b}, B\"ohm~\cite{Bohm1999}, and unpublished work of Bryant (see, for example,~\cite[Chapter 1]{Chow_etal}).  Here, we will present a sketch of the proof of the classification of these examples which emphasizes how one can view them as a one-parameter family of quasi-Einstein metrics parameterized by the dimensional parameter $m$.

Before listing the examples, we first recall that a BER flat SMMS $(M^n,g,v^m\dvol)$ with \charconstant\ $\mu$ satisfies
\begin{align*}
\Delta\phi - |\nabla\phi|^2 & = -m\mu e^{2\phi/m} & \text{if $m<\infty$ } \\
\Delta\phi - |\nabla\phi|^2 & = -\mu^\prime & \text{if $m=\infty$} .
\end{align*}
In~\cite{Case2009}, it is shown that if $g$ is complete, then either $\phi$ is constant and $\mu,\mu^\prime=0$, or $\mu,\mu^\prime>0$.  Thus the \charconstant\ of the examples we discuss here is necessarily positive.

As mentioned above, these manifolds have already been classified.  In the case $n=2$, these SMMS can be written very explicitly, which in particular clearly illustrates that these examples form a smooth one-parameter family with $m\in(1,\infty]$.

\begin{thm}[\cite{Berard-Bergery1982,Besse,Hamilton1986b}]
\label{thm:cigar}
Given $m\in(1,\infty]$, the unique rotationally symmetric quasi-Einstein SMMS on $\bR^2$ with quasi-Einstein constant $\lambda=0$ and \charconstant\ $\mu=\frac{4}{m-1}$ is
\begin{equation}
\label{eqn:cigar}
\left( \bR^2, dt^2\oplus\left(\frac{m-1}{2}v^\prime(t)\right)^2d\theta^2, v^m(t)\,\dvol \right),
\end{equation}
where $v(t)$ is the unique solution to the ODE
\begin{equation}
\label{eqn:cigar_ode}
\left(\frac{m-1}{2}v^\prime(t)\right)^2 = 1-v^{1-m}(t)
\end{equation}
with $v(0)=1$ and $v^\prime(t)\geq0$ for $t\geq 0$.
\end{thm}

The main point of the proof is that if one assumes that $(\bR^2,dt^2\oplus w^2\,d\theta^2,v^m\dvol)$ is a quasi-Einstein SMMS, it must be the case that $w=cv^\prime$ for some constant $c>0$.  The quasi-Einstein assumption then reduces to a second order ODE in $v$, and integrating this ODE once using the desired initial conditions yields~\eqref{eqn:cigar_ode}.

The point we wish to make here is that the normalizations we have taken in Theorem~\ref{thm:cigar} ensure that one has a well-defined limit as $m\to\infty$.  Indeed, recalling that $v=e^{-\frac{\phi}{m}}$, we see that
\[ \lim_{m\to\infty} \frac{m-1}{2}v^\prime(t) = -\frac{1}{2}\phi^\prime(t), \quad \lim_{m\to\infty} 1-v^{1-m}(t) = 1-e^{\phi(t)} . \]
In particular, the solution to~\eqref{eqn:cigar_ode} when $m=\infty$ with the specified initial conditions is $\phi=\log\sech^2t$, whence~\eqref{eqn:cigar} can be written
\[ \left( \bR^2, dt^2\oplus\tanh^2t\,d\theta^2, \cosh^2t\,\dvol, \infty\right), \]
which is the usual form of Hamilton's cigar soliton~\cite{Hamilton1986b}.  On the other hand, the quasi-Einstein SMMS~\eqref{eqn:cigar} when $2\leq m\in\bN$ are just the two-dimensional bases of the $(m+1)$-dimensional Riemannian Schwarzschild metric (cf.\ \cite{Besse}).

In the case $n\geq 3$, one does not have such an explicit form of the SMMS.  Nevertheless, the solutions have been completely classified by B\"ohm~\cite{Bohm1999} in the case $m\in(1,\infty)$ and by Bryant (see~\cite{Chow_etal}) in the case $m=\infty$.  By carefully checking the details of both proofs, one can also verify that these SMMS form a smooth one-parameter family in $m\in(1,\infty]$.

\begin{thm}
\label{thm:bryant}
Given $m\in(1,\infty]$, there exist non-constant functions $\phi,\psi\in C^\infty(M)$ such that the SMMS
\begin{equation}
\label{eqn:bryant_smms}
\left( \bR^n, dr^2\oplus\psi^2(r)d\theta^2, e^{-\phi}\dvol, m \right)
\end{equation}
is a nontrivial quasi-Einstein SMMS with quasi-Einstein constant zero.  Moreover, there is a unique such solution with $\phi(0)=1$ and \charconstant\ $\mu=\frac{1}{m-1}$.  With this choice,
\[ \psi^2(r) \sim \frac{n-2}{m+n-2}r^2 \]
for large $r$ when $m<\infty$ and
\[ \psi^2(r) \sim c(n)r \]
for large $r$ when $m=\infty$.
\end{thm}

\begin{proof}

The following proof is essentially due to B\"ohm~\cite{Bohm1999}, with coordinates changed in order to invoke the derivation of the Bryant solitons (cf.\ \cite{Chow_etal,Ivey_1994}).  First we observe that, when $m<\infty$, finding the desired BER-flat metric is equivalent to solving the equation
\begin{equation}
\label{eqn:qe}
\Ric - \frac{m}{v}\nabla^2 v = 0,
\end{equation}
where we have written $v=e^{-\phi/m}$ as usual.  In order for the metric and the density to be smooth at $r=0$, we require
\begin{equation}
\label{eqn:constraints}
\psi(0)=0, \qquad \psi^\prime(0) = 1, \qquad v^\prime(0)=0 .
\end{equation}
Additionally, \eqref{eqn:bianchi} yields the integrability condition
\begin{equation}
\label{eqn:constraint}
m\mu\frac{\psi^2}{v^2} = 2m(n-1)\psi\psi^\prime\frac{v^\prime}{v} + m(m-1)\psi^2\frac{(v^\prime)^2}{v^2} + (n-1)(n-2)[(\psi^\prime)^2-1] .
\end{equation}
Recalling that $v^\prime/v=-\phi^\prime/m$ and $m\mu_m\to\mu_\infty^\prime$, we see that these equations include the case $m=\infty$.

Now, \eqref{eqn:qe} and~\eqref{eqn:constraint} possess as symmetries translation of the coordinate $r$, simultaneous scaling of $v^2$ and $\mu$, and simultaneous scaling of $r$, $\psi$ and $\mu^{-1/2}$.  We thus introduce new coordinates which remove these symmetries, and for which, by~\eqref{eqn:constraint}, the solution is a curve on some compact submanifold of some $\bR^k$.  This is accomplished by the coordinates
\begin{align*}
X & = \sqrt{\frac{m+n-2}{(m-1)(n-2)}}\psi^\prime Y \\
Y & = \sqrt{\frac{(m-1)(n-1)(n-2)}{m}}\left((m-1)\psi v^{-1}v^\prime + (n-1)\psi^\prime\right)^{-1} \\
W & = \sqrt{\frac{m\mu}{(n-1)(n-2)}}\psi v^{-1} Y \\
dr & = \sqrt{\frac{m}{(m-1)(n-1)(n-2)}}\psi Y \; dt .
\end{align*}
The constants above are such that~\eqref{eqn:constraint} implies that the solution lies on the unit sphere $X^2+Y^2+W^2=1$.  One then checks that~\eqref{eqn:qe} becomes
\begin{subequations}
\label{eqn:ode2}
\begin{align}
\dot X & = X\left(X^2-\frac{1}{m-1}Y^2\right) + \frac{1}{m-1}\sqrt{\frac{m(m+n-2)}{n-1}}Y^2 - X \\
\dot Y & = Y\left(X^2-\frac{1}{m-1}Y^2\right) - \frac{1}{m-1}\sqrt{\frac{m(m+n-2)}{n-1}}XY + \frac{1}{m-1}Y \\
\dot W & = W\left(X^2-\frac{1}{m-1}Y^2\right) ,
\end{align}
\end{subequations}
where ``dot'' denotes differentiation with respect to $t$.  Moreover, the initial condition~\eqref{eqn:constraints} becomes
\[ (X,Y,W) \to I:=\left(\sqrt{\frac{m+n-2}{m(n-1)}},\sqrt{\frac{(m-1)(n-2)}{m(n-1)}},0\right) \qquad\qquad \mbox{ as } t\to-\infty . \]
As we require $v,\psi>0$, the desired solution is also constrained to lie in the first octant $X,Y,W\geq 0$.  A simple calculation then shows that, in addition to $I$, the only fixed point of~\eqref{eqn:ode2} on the unit sphere and in the first octant is
\[ K := \left(\sqrt{\frac{n-1}{m(m+n-2)}}, \sqrt{\frac{(m-1)(n-1)}{m(m+n-2)}}, \sqrt{\frac{m-1}{m+n-2}} \right) . \]
Indeed, Theorem~\ref{thm:bryant} is equivalent to the statement that there is a unique solution to~\eqref{eqn:ode2} which emerges from $I$ and tends to $K$, together with a study of its asymptotic behavior.  Note here that, in the case $m=\infty$, these are the coordinates normally used to derive the Bryant soliton~\cite{Chow_etal}, where one ignores the coordinate $W$ and wants to show that the solution $(X,Y)$ tends to $(0,0)$.

In order to prove such a solution converges to $K$, we introduce the Lyapunov function
\[ \kappa = W^{-\frac{2m}{m+n-1}} Y^{-\frac{2(n-1)}{m+n-1}} \left(1-\sqrt{\frac{n-1}{m(m+n-2)}}X\right)^2 \]
discovered by B\"ohm~\cite{Bohm1999}.  A straightforward computation using~\eqref{eqn:ode2} yields
\begin{align*}
\frac{\kappa^\prime}{\kappa} & = -\frac{2(n-1)}{(m-1)(m+n-1)}\left(1-\sqrt{\frac{m(m+n-2)}{n-1}}X\right)^2 \\
& \quad \times \left(1-\sqrt{\frac{n-1}{m(m+n-2)}}X\right)^{-1} ,
\end{align*}
and so $\kappa$ is indeed a Lyapunov function in the first octant $X,Y,W>0$.  Hence any solution in the first octant must converge to $K$.  Note also that, since $X^2+Y^2+W^2=1$, we can write $\kappa$ in terms of $X$ and $Y$ only.  In the case $m=\infty$, this implies that
\[ \kappa = W^{-2} = -\frac{1}{X^2+Y^2-1} . \]

In order to prove the uniqueness assertion of Theorem~\ref{thm:bryant}, we need to consider the critical point $I$ from which the solution emerges.  We first observe that the linearization of~\eqref{eqn:ode2} at $I$ is
\begin{align*}
\dot x & = \frac{2(m+n-2)-m(n-2)}{m(n-1)}x + \frac{2\sqrt{(m-1)(m+n-2)(n-2)}}{m(n-1)}y \\
\dot y & = \frac{(m-2)\sqrt{(m-1)(m+n-2)(n-2)}}{m(m-1)(n-1)}x + \frac{2(m+n-2)-2m(n-1)}{m(m-1)(n-1)}y \\
\dot w & = \frac{1}{n-1}w.
\end{align*}
Hence $I$ is a saddle, as the above system has eigenvalues
\[ \frac{1}{n-1}, \frac{2}{n-1} \mbox{ and } -\frac{n-2}{n-1} \]
whose corresponding eigenvectors are multiples of
\[ (0,0,1)\quad,\quad \left(2\sqrt{(m-1)(m+n-2)(n-1)},(m-2)(n-1),0\right), \]
and
\[ \left(\sqrt{(m-1)(n-2)},-\sqrt{m+n-2},0 \right) . \]
Since only the first and third eigenvectors are tangent to the unit sphere at $I$, there is only one relevant unstable direction, and hence there is a unique nonconstant solution which emerges from $I$.  Combined with the previous paragraph, we have shown that there is a unique solution $(X(t),Y(t),Z(t))$ in the first octant such that $(X,Y,Z)\to I$ as $t\to-\infty$.  Since fixing $\phi(0)$ and $\mu=\frac{1}{m-1}$ remove the symmetries of~\eqref{eqn:qe} and~\eqref{eqn:constraint}, this yields the claimed uniqueness.

Next, we consider the asymptotic behavior of the solutions.  Returning to our original coordinates, we see that when $m<\infty$, the critical point $K$ corresponds to
\[ (\psi^\prime,v^\prime) = \left(\sqrt{\frac{n-2}{m+n-2}},\sqrt{\frac{\mu}{m+n-2}}\right) . \]
In particular, we see that
\[ \psi(r)\sim\sqrt{\frac{n-2}{m+n-2}}r \]
for large $r$.  By the completeness of the system~\eqref{eqn:ode2} for $t$, this implies that, in the coordinates $(\psi,v,r)$, our solution is defined for all $r\geq 0$, and hence is asymptotically conical:
\[ g \sim dr^2 \oplus \frac{n-2}{m+n-2}r^2 d\theta_{n-1}^2 \mbox{ for } r\sim\infty . \]

To understand the behavior when $m=\infty$, we first observe that for all $1<m\leq\infty$, it holds that
\begin{equation}
\label{eqn:limiting}
\frac{X}{Y^2} \to \frac{1}{m-1}\sqrt{\frac{m(m+n-2)}{n-1}} .
\end{equation}
Next, when $m=\infty$, we have the relations
\begin{align}
\label{eqn:rel1}
\frac{d}{dY}(\log\psi+\log Y) & = \frac{X}{Y(X-\alpha)} \\
\label{eqn:rel2}
-\frac{dr}{dY} & = \frac{\psi}{\sqrt{(n-1)(n-2)}(X^2-\alpha X)},
\end{align}
where $\alpha=1/\sqrt{n-1}$.  Using the limiting behavior~\eqref{eqn:limiting}, the first relation implies that $\psi=O(Y^{-1})$, and the second implies that $r=O(Y^{-2})$.  Thus $\psi^2=O(r)$, as desired.

Inspecting the above proof carefully, we observe that all quantities are uniformly bounded in $m$, provided $m$ is uniformly bounded away from $1$.  Thus standard ODE theory implies that the quasi-Einstein metrics constructed by Theorem~\ref{thm:bryant} are continuous in $m\in(1,\infty]$.
\end{proof}

Finally, we remark that B\"ohm~\cite{Bohm1999} in fact proved that there is an $r$-parameter family of BER-flat metrics and an $(r+1)$-parameter family of metrics satisfying $\Ric_\phi^m=-g$ on the manifolds
\begin{equation}
\label{eqn:mfd}
\bR^n \times F_1 \times \cdots \times F_r
\end{equation}
for $F_i$ Einstein manifolds with positive scalar curvature, which are all rotationally symmetric multiply warped products with base $\bR^n$.  On the other hand, Dancer and Wang~\cite{DancerWang2008} have shown that there is an $r$-parameter family of steady gradient Ricci solitons and an $(r+1)$-parameter family of expanding gradient Ricci solitons which are rotationally symmetric warped product metrics on the manifolds~\eqref{eqn:mfd}.  It is thus natural to expect a similar statement to Theorem~\ref{thm:bryant} holds in this case, though the author has not checked these details.
\subsection{The L\"u-Page-Pope examples}
\label{sec:lpp}

In this section, we briefly discuss an interesting family of quasi-Einstein SMMS which are $S^2$-bundles over K\"ahler-Einstein manifolds.  These examples originate with Page's nonhomogeneous Einstein metric on $\bC P^2 \connectsum \overline{\bC P^2}$~\cite{Page1979}.  This construction was later rewritten in the language of cohomogeneity one metrics by B\'erard-Bergery~\cite{Berard-Bergery1982}, producing examples of $S^2$-bundles over K\"ahler-Einstein manifolds with positive scalar curvature.  These examples were soon rederived by Page and Pope~\cite{PagePope1987} using Page's original point of view.  They then appeared in the context of gradient Ricci solitons in the work of Koiso~\cite{Koiso1990}, H.-D.\ Cao~\cite{Cao1996}, and Chave and Valent~\cite{ChaveValent1996}.  Finally, L\"u, Page and Pope~\cite{LPP} completed the picture by constructing the warped product examples.

Like in the proof of Theorem~\ref{thm:lpp}, our main point is to illustrate how the examples of L\"u, Page and Pope~\cite{LPP} and the example of Koiso~\cite{Koiso1990}, Cao~\cite{Cao1996}, and Chave and Valent~\cite{ChaveValent1996} can be realized as a smooth one-parameter family of quasi-Einstein SMMS with $m\in(1,\infty]$.  Since these are examples of nontrivial compact quasi-Einstein SMMS, these examples are particularly relevant in motivating the precompactness theorem of~\cite{Case2010b}.

\begin{thm}
\label{thm:lpp}
Let $(M^{n-2},h)$ be a K\"ahler-Einstein manifold with K\"ahler form $\omega$ such that $\Ric=nh$, and suppose that the first Chern class $c_1=\frac{n}{2\pi}[\omega]$ is equal to $q\alpha^{\bR}$, where $\alpha\in H^2(M,\bZ)$ is an indivisible class, $q\in\bN$, and $\alpha^{\bR}$ is the real part of $\alpha$.  Then there are $q-1$ one-parameter families of quasi-Einstein SMMS with positive quasi-Einstein constant on $S^2$-bundles over $M$ which are parameterized by $m\in(1,\infty]$.
\end{thm}

\begin{remark}
It is well-known that $1\leq q\leq\frac{n}{2}$, with $q=\frac{n}{2}$ if and only if $M=\bC P^{\frac{n-2}{2}}$.
\end{remark}

\begin{proof}

The following sketch of the proof is an adaptation of the sketch found in~\cite{Besse}.  Further details can be found in the original references~\cite{Berard-Bergery1982,Cao1996,ChaveValent1996,Koiso1990,LPP}.

Let $s\in\bN$ be such that $1\leq s<q$.  Let $P(s)$ denote the principal $S^1$-bundle classified by $s\alpha$, and let $g(a,b)$ be the unique metric on $P(s)$ such that the projection $\pi\colon P(s)\to M$ is a Riemannian submersion onto $(M,b^2s)$ with totally geodesic fibers $S^1$ of length $2\pi a$.

To find the desired metrics, let $M(s)=[0,l]\times P(s)$ with the identification given by $\pi$ on the boundary.  Let $t$ be the standard coordinate on $[0,l]$ and consider the metric $\og=dt^2\oplus g(f(t),h(t))$ on $M(s)$ for $f,h$ positive functions on $(0,l)$.  For this metric to be smooth, we require that
\[ f(0)=0=f(l), \quad f^\prime(0)=1=-f^\prime(l), \quad h>0 \mbox{ on } [0,l], \quad \mbox{and } h^\prime(0)=0=h^\prime(l). \]
Fix $m$, let $v=e^{-\phi/m}$ , and assume $\Ric_\phi^m=\lambda g$.  If $m<\infty$, this implies that
\begin{align*}
\lambda & = -\frac{f^{\prime\prime}}{f} - (n-2)\frac{h^{\prime\prime}}{h} - m\frac{v^{\prime\prime}}{v} \\
\lambda & = - \frac{f^{\prime\prime}}{f} - (n-2)\frac{f^\prime h^\prime}{fh} - m\frac{f^\prime v^\prime}{fv} + \frac{(n-2)n^2s^2f^2}{4q^2 h^4} \\
\lambda & = -\frac{h^{\prime\prime}}{h} - \frac{f^\prime h^\prime}{fh} - m\frac{h^\prime v^\prime}{hv} - (n-3)\left(\frac{h^\prime}{h}\right)^2 + \frac{n}{h^2} - \frac{n^2s^2 f^2}{2q^2 h^4} ,
\end{align*}
and if $m=\infty$, simply use the facts that $-m\frac{v^\prime}{v}\to \phi^\prime$, $-m\frac{v^{\prime\prime}}{v}\to \phi^{\prime\prime}$ as $m\to\infty$ to derive the equivalent system of ODEs.  From~\eqref{eqn:bianchi}, we have the integrability condition
\begin{align*}
\mu & = \lambda v^2 + v v^{\prime\prime} + \frac{vf^\prime v^\prime}{f} + (n-2)\frac{vh^\prime v^\prime}{h} + (m-1)(v^\prime)^2 & \mbox{ if } m& <\infty, \\
-\mu^\prime+n\lambda & = 2\lambda\phi + \phi^{\prime\prime} + \frac{f^\prime}{f}\phi^\prime + (n-2)\frac{h^\prime}{h}\phi^\prime - (\phi^\prime)^2 & \mbox{ if } m& =\infty .
\end{align*}
Of course, by Proposition~\ref{prop:bianchi_limit}, the latter equation is the limit of the former as $m\to\infty$, as seen by multiplying by $-mv^{-2}$.  One can solve this system by setting $\mu=m-1$ if $m<\infty$ and $\mu=0$ if $m=\infty$, and moreover, can show that the solution will be complete.
\end{proof}

The derivation of these metrics is presented for gradient Ricci solitons in the more general setting of cohomogeneity one metrics by Dancer and Wang~\cite{DancerWang2008}, where they consider even more general bundles over the base manifold $M$.  Using these more general bundles in the same way as above, it is reasonable to expect that one can find even more families of quasi-Einstein metrics (cf.\ \cite{Chen2009,WangWang1998}).

\begin{remark}
It is interesting to note that the examples constructed by L\"u, Page and Pope~\cite{LPP} are not K\"ahler, while the limits constructed by Cao~\cite{Cao1996}, Koiso~\cite{Koiso1990} and Chave and Valent~\cite{ChaveValent1996} are K\"ahler.  This should be compared with a theorem of Shu, Wei and the author~\cite{Case_Shu_Wei}, which states that there are no nontrivial quasi-Einstein SMMS with $m\in(0,\infty)$ and positive quasi-Einstein constant for which the underlying metric is K\"ahler.  In light of the perspective of the present article, it would be interesting to know if the L\"u-Page-Pope examples are conformally K\"ahler.
\end{remark}
\subsection{Constructing New Examples via Products}
\label{sec:products}

We conclude this article by discussing a two product constructions for quasi-Einstein SMMS, generalizing the basic fact that the trivial product of two Einstein metrics with the same Einstein constant is itself Einstein.  In particular, this construction allows us to produce examples of quasi-Einstein SMMS with $m=1$ and nonzero \charconstant\ (cf.\ Remark~\ref{rk:m=1_trouble}).

From the perspective of the \auxiliarymanifold, our two product constructions arise by taking a product either with $M$ or $F$.  More precisely, given a SMMS $(M^n,g,v^m\dvol)$ with $2\leq m\in\bN$, there are two interesting ways to form products.  First, we can take a product with the entire \auxiliarymanifold; i.e.\ with $(F^m,h)$ denoting the fiber and $(N^k,g_2)$ denoting some additional Riemannian manifold, we can consider the warped product manifold
\[ \left( M^n\times N^k\times F^m, g\oplus g_2\oplus v^2h\right) . \]
Second, we can instead think of regarding the fiber as a product $(N^k\times F^{m-k},g_2\oplus h)$, and regard the \auxiliarymanifold\ as
\[ \left( M^n\times N^k\times F^{m-k}, g\oplus v^2g_2\oplus v^2h\right) . \]
In both cases, we then think of collapsing the fibers $F$ as in Example~\ref{ex:cc} to yield the desired product construction.  We formalize this to an intrinsic statement on SMMS as follows.

\begin{defn}
Let $(M^n,g,v^m\dvol_g)$ be a SMMS and let $(N^k,h)$ be a Riemannian manifold.  The \emph{flat product of $M$ and $N$} is the SMMS
\[ \left( M^n\times N^k, \og:=g\oplus h, v^m\dvol_{\og}\right) . \]
\end{defn}

\begin{defn}
Let $(M^n,g,v^m\dvol_g)$ be a SMMS and let $(F^k,h)$ be a Riemannian manifold.  The \emph{warped product of $M$ and $N$} is the SMMS
\[ \left( M^n\times N^k, \og:=g\oplus v^2h, v^{m-k}\dvol_{\og} \right) . \]
\end{defn}

As stated above, one immediately notices that the flat product leaves the dimensional parameter $m$ unchanged, while the warped product decreases the dimensional parameter $k$ by the dimension of the fiber.  In fact, if $m=k\in\bN$, taking the warped product as above is simply constructing the \auxiliarymanifold.  The key points here are that both definitions do not impose any constraints on $m$ (though observe that the two products coincide when $\lv m\rv=\infty$), and that they allow one to construct new quasi-Einstein SMMS by taking products of certain types of SMMS and Riemannian manifolds:

\begin{prop}
Let $(M^n,g,v^m\dvol)$ be a quasi-Einstein SMMS with quasi-Einstein constant $\lambda$ and \charconstant\ $\mu$.  If $(F^k,h)$ is an Einstein manifold with $\Ric(h)=\lambda h$, then the flat product of $M$ and $F$ is a quasi-Einstein SMMS with quasi-Einstein constant $\lambda$ and \charconstant\ $\mu$.
\end{prop}

\begin{prop}
Let $(M^n,g,v^m\dvol)$ be a quasi-Einstein SMMS with quasi-Einstein constant $\lambda$ and \charconstant\ $\mu$.  If $(F^k,h)$ is an Einstein manifold with $\Ric(h)=\mu h$, then the warped product of $M$ and $F$ is a quasi-Einstein SMMS with quasi-Einstein constant $\lambda$ and \charconstant\ $\mu$.
\end{prop}

In both cases, the proof is a straightforward computation using the standard formulas for the Ricci curvature and Hessian of a warped product (see~\cite{Besse}).  As a consequence of these results, we can construct the aforementioned quasi-Einstein SMMS with $m=1$ and nonzero \charconstant:

\begin{cor}
Let $(M^n,g,v^m\dvol)$ be a nontrivial compact quasi-Einstein SMMS with $3\leq m\in\bN$ and \charconstant\ $\mu=m-1$, and let $(S^{m-1},h)$ be the standard sphere with sectional curvature one.  Then the warped product of $M$ and $S^{m-1}$ is a nontrivial compact quasi-Einstein SMMS with dimensional parameter $m=1$ and \charconstant\ $\mu\not=0$.
\end{cor}

\begin{remark}

The examples of L\"u, Page and Pope~\cite{LPP} give SMMS which meet the hypotheses of the corollary.
\end{remark}             

\bibliographystyle{abbrv}
\bibliography{../bib}
\end{document}